\documentclass[a4paper,10pt,twoside]{article}
\usepackage{pst-fill,pst-grad}
\usepackage{textcomp}
\usepackage[english]{babel}
\usepackage[utf8x]{inputenc}
\usepackage{graphicx}
\usepackage{amsmath}
\usepackage{float}
\usepackage{fancyhdr}
\usepackage[matrix,arrow,curve]{xy}
\usepackage{pstricks} 
\usepackage{amsmath,amsfonts,verbatim,afterpage,theorem,euscript,mathrsfs,amssymb}
\usepackage{amsfonts}
\usepackage{amssymb}
\usepackage{array}
\usepackage{dsfont}
\usepackage{hyperref}
\newtheorem{Definition}{Definition}[section]
\newtheorem{Proposition}{Proposition}[section]
\newtheorem{Lemme}{Lemma}[section]
\newtheorem{Theoreme}{Theorem}
\newtheorem{Corollaire}{Corollary}[section]
\newtheorem{Remarque}{Remark}[section]
\title{\bf Remarks on a fractional diffusion transport equation\\ with applications to the\\ dissipative quasi-geostrophic equation}
\author{Diego Chamorro\footnote{Laboratoire d'Analyse et de Probabilit\'es, Universit\'e d'Evry Val d'Essonne, 23 Boulevard de France, 91037 Evry Cedex - France, diego.chamorro@univ-evry.fr}}
\begin{document}
\maketitle
\begin{scriptsize}
\abstract{In this article we study  Hölder regularity $\mathcal{C}^\gamma$ for solutions of a transport equation based in the dissipative quasi-geostrophic equation. Adapting an idea of A. Kiselev and F. Nazarov presented in \cite{KN}, we will use the molecular characterization of local Hardy spaces $h^\sigma$ in order to obtain information on Hölder regularity of such solutions. This will be done by following the evolution of molecules in a backward equation. We will also study global existence, Besov regularity for weak solutions and a maximum principle  and we will apply these results to the critical dissipative quasi-geostrophic equation.}\\[3mm]
\textbf{Keywords: quasi-geostrophic equation, Hölder regularity, Hardy spaces.}
\end{scriptsize}
\section{Introduction}
The dissipative quasi-geostrophic equation has been studied by many authors, not only because of its own mathematical importance, but also as a 2D model in geophysical fluid dynamics and because of its close relationship with other equations arising in fluid dynamics. See \cite{CW}, \cite{Marchand1} and the references given there for more details. This equation has the following form:
\begin{equation*}
(QG)_{\alpha}\quad \left\lbrace
\begin{array}{rl}
\partial_t \theta(x,t)=& \nabla\cdot(u\,\theta)(x,t)-\Lambda^{2\alpha}\theta(x,t)\\[5mm]
\theta(x,0)=& \theta_0(x)
\end{array}
\right.
\end{equation*}
for $0<\alpha\leq 1$ and $t\in [0, T]$. Here $\theta$ is a real-valued function and the velocity $u$ is defined by means of Riesz transforms in the following way:
$$u=(-R_2\theta, R_1\theta).$$
Recall that Riesz Transforms $R_j$ are given by $\widehat{R_j \theta}(\xi)=-\frac{i\xi_j}{|\xi|}\widehat{\theta}(\xi)$ for $j=1,2$ and that
$\Lambda^{2\alpha}=(-\Delta)^{\alpha}$ is the Laplacian's fractional power defined by the formula
$$\widehat{\Lambda^{2\alpha} \theta}(\xi)=|\xi|^{2\alpha}\widehat{\theta}(\xi)$$
where $\widehat{\theta}$ denotes the Fourier transform of $\theta$.\\

It is classical to consider three cases in the analysis of the dissipative quasi-geostrophic equation following the values of the diffusion parameter $\alpha$. The case $1/2<\alpha$ is called \textit{sub-critical} since the diffusion factor is stronger than the nonlinearity. In this case, weak solutions were constructed by S. Resnick in \cite{Resnick} and P. Constantin \& J. Wu showed in \cite{CW1} that smooth initial data gives a smooth global solution. 

The \textit{critical} case is given when $\alpha=1/2$. Here, P. Constantin, D. C\'ordoba \& J. Wu studied in \cite{CCW} global existence in Sobolev spaces, while global well-posedness in Besov spaces has been treated by H. Abidi \& T. Hmidi in \cite{Hm}. Also in this case, and more recently, A. Kiselev, F. Nazarov \& A. Volberg showed in \cite{KNV} that any regular periodic data generates a unique $\mathcal{C}^\infty$ solution. 

Finally, the case when $0<\alpha<1/2$ is called \textit{super-critical} partially because it is harder to work with than the two other cases, but mostly because the diffusion term is weaker than the nonlinear term. In this last case, weak solutions for initial data in $L^p$ or in $\dot{H}^{-1/2}$ were studied by F. Marchand in \cite{Marchand}.\\

In this article, following L. Caffarelli \& A. Vasseur in \cite{Caffarelli}, we will study a special version of the dissipative quasi-geostrophic equation  $(QG)_\alpha$. The idea is to replace the Riesz Transform-based velocity $u$ by a new velocity $v$ to obtain the following $n$-dimensional fractional diffusion transport equation for $0<\alpha\leq 1$:
\begin{equation}\label{QG}
(T)_\alpha \quad \left\lbrace
\begin{array}{rl}
\partial_t \theta(x,t)=& \nabla\cdot(v\,\theta)(x,t)-\Lambda^{2\alpha}\theta(x,t)\\[5mm]
\theta(x,0)=& \theta_0(x)\\[5mm]
div(v)=&0.
\end{array}
\right.
\end{equation}
In (\ref{QG}), the functions $\theta$ and $v$ are such that $\theta:\mathbb{R}^n\times [0,T]\longrightarrow \mathbb{R}$ and $v:\mathbb{R}^n\times [0,T]\longrightarrow \mathbb{R}^n$. Remark that the velocity $v$ is now a given data for the problem and we will always assume that $v$ is divergence free and belongs to $L^{\infty}([0,T];bmo(\mathbb{R}^n))$. \\

We fix once and for all the parameter $\alpha=1/2$ in order to study the \textbf{critical} case. The main theorem presented in this article studies the regularity of the solutions of the fractional diffusion transport equation $(T)_{1/2}$:
\begin{Theoreme}[Hölder regularity]\label{Theo1}
Let $\theta_0$ be a function such that $\theta_0\in L^{\infty}(\mathbb{R}^n)$ and $T_0>0$ a small positive time. If $\theta(x,t)$ is a solution for the equation $(T)_{1/2}$, then for all time $T_0<t<T$,  we have that $\theta(\cdot,t)$ belongs to the H\"older space $\mathcal{C}^{\gamma}(\mathbb{R}^n)$ with $0<\gamma<1$. 
\end{Theoreme}
This result says that after a time $T_0$ there is a small smoothing effect so the dissipation given by the fractional Laplacian is stronger than the drift term in equation $(T)_{1/2}$. Thus, since the velocity $v$ in equation $(T)_{1/2}$ belongs to $L^{\infty}([0,T];bmo(\mathbb{R}^n))$, it would be quite simple to adapt this result to the $(QG)_{1/2}$ equation. See section \ref{SecAppli} for details.\\

Let us say a few words about the proof of this theorem. A classical result of harmonic analysis states that Hölder spaces $\mathcal{C}^\gamma(\mathbb{R}^n)$ can be paired with local Hardy spaces $h^\sigma(\mathbb{R}^n)$. Therefore, if we prove that the duality bracket 
\begin{equation}\label{DualQuantity}
\langle  \theta(\cdot,t),\psi_0\rangle=\int_{\mathbb{R}^n} \theta(x,t)\psi_0(x)dx
\end{equation}
is bounded for every $\psi_0\in h^\sigma(\mathbb{R}^n)$, with $t\in ]0,T[$, we obtain that $\theta(\cdot,t)\in \mathcal{C}^{\gamma}(\mathbb{R}^n)$. 
One of the main features of Hardy spaces is that they admit a characterization by \textit{molecules} (see definition \ref{DefMolecules} below), which are rather simple functions, and this allows us to study the quantity (\ref{DualQuantity}) only for such molecules.\\

This dual approach was originally given in the torus $\mathbb{T}^n$ by A. Kiselev \& F. Nazarov in \cite{KN} with a very special family of test functions. Thus, the main novelty of this paper besides the generalization to $\mathbb{R}^n$ is the use of molecular Hardy spaces.\\

Broadly speaking and following \cite{Stein2} p. 130, a \emph{molecule} is a function $\psi$ so that 
\begin{enumerate} 
\item[(i)] $\displaystyle{\int_{\mathbb{R}^n}}\psi(x)dx=0$\\
\item[(ii)] $|\psi(x)|\leq  r^{-n/\sigma}\min\{1; r^{\beta n}/ |x-x_0|^{\beta n}\}$,
\end{enumerate}
with $\beta \sigma >1$ and $x_0\in \mathbb{R}^n$, where the parameter $r\in ]0, +\infty[$ stands for the \textit{size} of the molecule $\psi$.\\ 

Since we are going to work with local Hardy spaces, we will introduce a size treshold in order to distinguish \textit{small} molecules from \textit{big} ones in the following way:
\begin{Definition}[$r$-molecules]\label{DefMolecules} Set $\frac{n}{n+1}<\sigma<1$, define $\gamma=n(\frac{1}{\sigma}-1)$ and fix a real number $\omega$ such that $0<\gamma<\omega<1$. An integrable function $\psi$ is an $r$-molecule if we have
\begin{enumerate}
\item[$\bullet$]\underline{Small molecules $(0<r<1)$:}
\begin{eqnarray}
& & \int_{\mathbb{R}^n} |\psi(x)||x-x_0|^{\omega}dx \leq  r^{\omega-\gamma}\mbox{, for } x_0\in \mathbb{R}^n\label{Hipo1} \;\qquad\qquad\qquad\mbox{(concentration condition)} \\[5mm]
& &\|\psi\|_{L^\infty}  \leq  \frac{1}{r^{n+\gamma}}\label{Hipo2} \qquad\qquad\qquad\qquad\qquad\qquad\qquad\qquad\qquad\, \mbox{(height condition)}  \\[5mm]
& &\int_{\mathbb{R}^n} \psi(x)dx=0\label{Hipo3}\qquad\qquad\qquad\qquad\qquad\qquad\qquad\qquad\qquad \mbox{(moment condition)} 
\end{eqnarray}
\item[$\bullet$] \underline{Big molecules $(1\leq r<+\infty)$:}\\

In this case we only require conditions (\ref{Hipo1}) and (\ref{Hipo2}) for the $r$-molecule $\psi$ while the moment condition (\ref{Hipo3}) is dropped.
\end{enumerate}
\end{Definition}
It is interesting to compare this definition of molecules to the one used in \cite{KN}. In our molecules the parameter $\gamma$ reflects explicitly the relationship between Hardy and H\"older spaces (see Theorem \ref{TheoHHD} below).  However, the most important fact relies in the parameter $\omega$ which gives us the additional flexibility that will be crucial in the following calculations.
\begin{Remarque}\label{Remark2}
\emph{
\begin{itemize}
\item[1)] Note that the point $x_0\in \mathbb{R}^n$ can be considered as the ``center'' of the molecule.
\item[2)] Conditions (\ref{Hipo1}) and (\ref{Hipo2}) are an easy consequence of condition (ii) and they both imply the estimate $\|\psi\|_{L^1}\leq C\, r^{-\gamma}$, thus every $r$-molecule belongs to $L^p(\mathbb{R}^n)$ for $1<p<+\infty$. We stress here that the previous $L^1$ estimate is just a corollary of (\ref{Hipo1}) and (\ref{Hipo2}). This inequality is not a part of the molecule's definition. 
\end{itemize}
}
\end{Remarque}
The main interest for using molecules relies in the possibility of \textit{transfering} the regularity problem to the evolution of such molecules:
\begin{Proposition}[Transfer property]\label{Transfert} Let $\psi(x,s)$ be a solution of the backward problem
\begin{equation}
\left\lbrace
\begin{array}{rl}
\partial_s \psi(x,s)=& -\nabla\cdot [v(x,t-s)\psi(x,s)]-\Lambda\psi(x,s)\label{Evolution1}\\[5mm]
\psi(x,0)=& \psi_0(x)\in L^1\cap L^\infty(\mathbb{R}^n)\\[5mm]
div(v)=&0 \quad \mbox{and }\; v\in L^{\infty}([0,T];bmo(\mathbb{R}^n))
\end{array}
\right.
\end{equation}
If $\theta(x,t)$ is a solution of (\ref{QG}) with $\theta_0\in L^\infty(\mathbb{R}^n)$ then we have the identity
\begin{equation*}
\int_{\mathbb{R}^n}\theta(x,t)\psi(x,0)dx=\int_{\mathbb{R}^n}\theta(x,0)\psi(x,t)dx.
\end{equation*}
\end{Proposition}
\textit{\textbf{Proof.}}
We first consider the expression
$$\partial_s\int_{\mathbb{R}^n}\theta(x,t-s)\psi(x,s)dx=\int_{\mathbb{R}^n}-\partial_s\theta(x,t-s)\psi(x,s)+\partial_s\psi(x,s)\theta(x,t-s)dx.$$
Using equations (\ref{QG}) and (\ref{Evolution1}) we obtain
\begin{eqnarray*}
\partial_s\int_{\mathbb{R}^n}\theta(x,t-s)\psi(x,s)dx&=&\int_{\mathbb{R}^n}- \nabla\cdot\left[(v(x,t-s)\theta(x,t-s)\right]\psi(x,s)+\Lambda\theta(x,t-s)\psi(x,s) \\[5mm]
&-& \nabla\cdot\left[(v(x,t-s) \psi(x,s))\right]\theta(x,t-s)-\Lambda\psi(x,s) \theta(x,t-s)dx.
\end{eqnarray*}
Now, using the fact that $v$ is divergence free we have that expression above is equal to zero, so the quantity
$$\int_{\mathbb{R}^n}\theta(x,t-s)\psi(x,s)dx$$
remains constant in time. We only have to set $s=0$ and $s=t$ to conclude. \hfill $\blacksquare$\\

This proposition says, that in order to control $\langle \theta(\cdot, t),\psi_0\rangle$, it is enough (and much simpler) to study the bracket $\langle \theta_0,\psi(\cdot, t)\rangle$. Let us explain in which sense this transfer property is useful: in the bracket $\langle \theta_0,\psi(\cdot, t)\rangle$ we have much more informations than in the bracket (\ref{DualQuantity}) since the initial data $\psi_0$ is a molecule which satisfies conditions (\ref{Hipo1})-(\ref{Hipo3}).\\

\textbf{Proof of the Theorem \ref{Theo1}.} 
Once we have the transfer property proven above, the proof of the Theorem \ref{Theo1} is quite direct and it reduces to a $L^1$ estimate for molecules. Indeed, assume that for \textit{all} molecular initial data $\psi_0$ we have a $L^1$ control for $\psi(\cdot,t)$ a solution of (\ref{Evolution1}), then the Theorem \ref{Theo1} follows easily: applying Proposition \ref{Transfert} with the fact that $\theta_0\in L^{\infty}(\mathbb{R}^n)$ we have 
\begin{equation}\label{DualQuantity1}
|\langle \theta(\cdot, t), \psi_0\rangle|=\left|\int_{\mathbb{R}^n}\theta(x,t)\psi_0(x)dx\right|=\left|\int_{\mathbb{R}^n}\theta(x,0)\psi(x,t)dx\right|\leq \|\theta_0\|_{L^\infty}\|\psi(\cdot,t)\|_{L^1}<+\infty.
\end{equation}
From this, we obtain that $\theta(\cdot, t)$ belongs to the Hölder space $\mathcal{C}^\gamma(\mathbb{R}^n)$.\\

Now we need to study the control of the $L^1$ norm of $\psi(\cdot, t)$ and we divide our proof in two steps following the molecule's size. For the initial big molecules, \textit{i.e.} if $r\geq 1$, the needed control is straightforward: apply the maximum principle (\ref{PrincipeMax2}) below and the remark \ref{Remark2}-2) to obtain
$$\|\theta_0\|_{L^\infty}\|\psi(\cdot,t)\|_{L^1}\leq \|\theta_0\|_{L^\infty}\|\psi_0\|_{L^1}\leq C \frac{1}{r^\gamma} \|\theta_0\|_{L^\infty},$$
but, since $r\geq 1$, we have that $|\langle \theta(\cdot, t), \psi_0\rangle|<+\infty$ for all \textit{big} molecules.\\

In order to finish the proof of the theorem, it  only remains to treat the $L^1$ control for \textit{small} molecules. This is the most complex part of the proof and we will present it in section \ref{SectL1control} where we will prove the next theorem:

\begin{Theoreme}\label{TheoL1estimate}
For all small initial molecular data $\psi_0$, there exists an small time $T_0>0$ such that:
$$\|\psi(\cdot, t)\|_{L^1}\leq C T_0^{-\gamma}\qquad \mbox{for all } T_0<t<T \qquad (0<\gamma<\omega<1).$$
\end{Theoreme}
Taking for granted this theorem, we obtain a good control over the quantity $\|\psi(\cdot, t)\|_{L^1}$  for all $0<r<1$. Finally, getting back to (\ref{DualQuantity1}) we obtain that $|\langle \theta(\cdot, t), \psi_0\rangle|$ is always bounded for $T_0<t<T$ and for any molecule $\psi_0$: we have proven by a duality argument the Theorem \ref{Theo1}.\hfill $\blacksquare$\\

Let us recall now that for (smooth) solutions of $(QG)_{1/2}$ and $(T)_{1/2}$ above we can use the remarkable property of \textit{maximum principle} mentioned before. This was proven in \cite{Cordoba} and in \cite{Marchand} and it gives us the following inequalities. 
\begin{eqnarray}
\|\theta(\cdot, t)\|_{L^p} +p\int_0^t\int_{\mathbb{R}^n}|\theta(x,s)|^{p-2}\theta(x,s) \Lambda \theta(x,s)dxds&\leq &\|\theta_0\|_{L^p}\qquad (2\leq p<+\infty)\label{PrincipeMax1}\\[3mm]
\mbox{or more generally }\qquad\qquad \|\theta(\cdot, t)\|_{L^p}&\leq& \|\theta_0\|_{L^p} \qquad (1\leq p\leq+\infty)\label{PrincipeMax2}
\end{eqnarray}
These estimates are extremely useful and they are the starting point of several works. Indeed, the study of inequality (\ref{PrincipeMax1}) helps us incidentally to solve a question pointed out by F. Marchand in \cite{Marchand} concerning weak solution's global regularity:
\begin{Theoreme}[Weak solution's regularity]\label{TheoRegBesov}
Let $2\leq p <+\infty$. If $\theta_0\in L^{p}(\mathbb{R}^n)$ is an initial data for $(QG)_{1/2}$ or $(T)_{1/2}$ equations, then the associated weak solution  $\theta(x,t)$ belongs to $L^\infty([0,T]; L^p(\mathbb{R}^n))\cap L^p([0,T]; \dot{B}^{1/p,p}_p(\mathbb{R}^n))$.
\end{Theoreme}
Theorem \ref{Theo1}, Theorem \ref{TheoRegBesov} and the $L^1$ control for small molecules are the core of the paper, however, for the sake of completness, we will prove some other interesting results concerning the equation (\ref{QG}).\\

The plan of the article is the following: in the section \ref{HardyHolderbmo} we recall some facts concerning the molecular characterization of local Hardy spaces and some other facts about Hölder and $bmo$ spaces. In section \ref{SectL1control} we study the $L^1$-norm control for molecules and in section \ref{SecExiUnic} we study existence and uniqueness of solutions with initial data in $L^p$ with $2\leq p<+\infty$ and we prove the Theorem \ref{TheoRegBesov}. Section \ref{Sect_PrincipeMax} is devoted to a \textit{positivity principle} that will be useful in our proofs and section \ref{SecLinfty} studies existence of solution with $\theta_0\in L^\infty$. Finally, section \ref{SecAppli} applies these results to the 2D-quasi-geostrophic equation $(QG)_{1/2}$.
\section{Molecular Hardy spaces, H\"older spaces and $bmo$}\label{HardyHolderbmo}

Hardy spaces have several equivalent characterizations (see \cite{Coifmann}, \cite{Gold} and \cite{Stein2} for a detailed treatment). In this paper we are interested mainly in the molecular approach that defines \textit{local} Hardy spaces $h^\sigma$ with $0<\sigma<1$: 
\begin{Definition}[Local Hardy spaces $h^\sigma$] Let $0<\sigma<1$. The local Hardy space $h^{\sigma}(\mathbb{R}^n)$ is the set of distributions $f$ that admits the following molecular decomposition:
\begin{equation}\label{MolDecomp}
f=\sum_{j\in \mathbb{N}}\lambda_j \psi_j
\end{equation}
where $(\lambda_j)_{j\in \mathbb{N}}$ is a sequence of complex numbers such that $\sum_{j\in \mathbb{N}}|\lambda_j|^\sigma<+\infty$ and $(\psi_j)_{j\in \mathbb{N}}$ is a family of $r$-molecules in the sense of the Definition \ref{DefMolecules} above. The $h^\sigma$-norm\footnote{it is not actually a \textit{norm} since $0<\sigma<1$. More details can be found in \cite{Gold} and \cite{Stein2}.} is then fixed by the formula 
$$\|f\|_{h^\sigma}=\inf\left\{\left(\sum_{j\in \mathbb{N}}|\lambda_j|^\sigma\right)^{1/\sigma}:\; f=\sum_{j\in \mathbb{N}}\lambda_j \psi_j \right\}$$ where the infimum runs over all possible decompositions (\ref{MolDecomp}).
\end{Definition}
Local Hardy spaces have many remarquable properties and we will only stress here, before passing to duality results concerning $h^\sigma$ spaces, the fact that Schwartz class $\mathcal{S}(\mathbb{R}^n)$ is dense in $h^{\sigma}(\mathbb{R}^n)$. For further details see \cite{Grafakos}, \cite{Stein2}, \cite{Gold} and \cite{CW}.\\

Now, let us take a closer look at the dual space of local Hardy spaces. In \cite{Gold} D. Goldberg proved the next important theorem:
\begin{Theoreme}[Hardy-Hölder duality]\label{TheoHHD} Let $\frac{n}{n+1}<\sigma<1$ and fix $\gamma=n(\frac{1}{\sigma}-1)$. Then the dual of local Hardy space $h^{\sigma}(\mathbb{R}^n)$ is the Hölder space $\mathcal{C}^\gamma(\mathbb{R}^n)$ fixed by the norm
$$\|f\|_{\mathcal{C}^\gamma}=\|f\|_{L^\infty}+\underset{x\neq y}{\sup}\frac{|f(x)-f(y)|}{|x-y|^\gamma}.$$
\end{Theoreme}
This result allows us to study the Hölder regularity of functions in terms of Hardy spaces and it will be applied to solutions of the $n$-dimensional fractional diffusion transport equation (\ref{QG}).
\begin{Remarque}\label{Remark3}\emph{Since $0<\sigma<1$, we have $\sum_{j\in \mathbb{N}}|\lambda_j|\leq\left(\sum_{j\in \mathbb{N}}|\lambda_j|^\sigma\right)^{1/\sigma}$ thus for testing Hölder continuity of a function $f$ it is enough to study the quantities $\langle f,\psi_j\rangle$ where $\psi_j$ is an $r$-molecule.}
\end{Remarque}
We finish this section by recalling some useful facts about the $bmo$ space used to characterize velocity $v$. This space is defined as locally integrable functions $f$ such that
$$\underset{|B|\leq 1}{\sup}\frac{1}{|B|}\int_{B}|f(x)-f_B|dx<M \qquad\mbox{ and }\qquad \underset{|B|>1}{\sup}\frac{1}{|B|}\int_{B}|f(x)|dx<M\qquad \mbox{for a constant } M;$$
 where we noted $B(R)$ a ball of radius $R>0$ and $f_B=\frac{1}{|B|}\displaystyle{\int_{B(R)}}f(x)dx$. The norm $\|\cdot\|_{bmo}$ is then fixed as the smallest constant $M$ satisfying these two conditions. We will use the next properties for a function belonging to $bmo$:
\begin{Proposition}\label{PropoBMO1} Let $f\in bmo$, then
\begin{enumerate}
\item[1)] for all $1<p<+\infty$, $f$ is locally in $L^p$ and $\frac{1}{|B|}\displaystyle{\int_{B}|f(x)-f_B|^pdx}\leq C \|f\|_{bmo}^p$
\item[2)] for all $k\in \mathbb{N}$, we have $|f_{2^k B}-f_B|\leq Ck \|f\|_{bmo}$ where $2^kB(R)=B(2^k R)$.
\end{enumerate}
\end{Proposition}
\begin{Proposition}\label{TheoApproxbmo}
Let $f$ be a function in $bmo(\mathbb{R}^n)$. For $k\in \mathbb{N}$, define $f_k$ by
\begin{equation}\label{Forbmoaprox}
f_k(x)=\left\lbrace
\begin{array}{rll}
-k & \mbox{if} & f(x)\leq -k \\[3mm]
f(x) & \mbox{if} & -k\leq f(x)\leq k \\[3mm]
k & \mbox{if} & k\leq f(x).
\end{array}
\right.
\end{equation}
Then $(f_k)_{k\in \mathbb{N}}$ converges weakly to $f$ in $bmo(\mathbb{R}^n)$.
\end{Proposition}
For a proof of these results and more details on Hardy, Hölder and $bmo$ spaces see \cite{Coifmann}, \cite{Gold}, \cite{PGLR}, \cite{Grafakos} and \cite{Stein2}.
\section{$L^1$ control for small molecules: proof of the Theorem \ref{TheoL1estimate}}\label{SectL1control}
As said in the introduction, we need to construct a suitable control in time for the $L^1$-norm of the solutions $\psi(\cdot,t)$ of the backward problem (\ref{Evolution1}) where the inital data $\psi_0$ is a \textit{small} $r$-molecule. This will be achieved by iteration in two different steps. The first step explains the molecules' deformation after a very small time $s_0>0$, which is related to the size $r$ by the bounds $0<s_0\leq \epsilon r$ with $\epsilon$ a small constant. In order to obtain a control of the $L^1$ norm for larger times we need to perform a second step which takes as a starting point the results of the first step and gives us the deformation for another small time $s_1$, which is also related to the original size $r$. Once this is achieved it is enough to iterate the second step as many times as necessary to get rid of the dependence of the times $s_0, s_1,...$ from the molecule's size. Proceeding this way we obtain the $L^1$ control needed for all time $T_0<t<T$.
\subsection{Small time molecule's evolution: First step}\label{SecEvolMol1}
The following theorem shows how the molecular properties are deformed with the evolution for a small time $s_0$.
\begin{Theoreme}\label{SmallGeneralisacion} Set $\sigma$, $\gamma$ and $\omega$ three real numbers such that $\frac{n}{n+1}<\sigma<1$, $\gamma=n(\frac{1}{\sigma}-1)$ and $0<\gamma<\omega<1$. Let $\psi(x,s_0)$ be a solution of the problem
\begin{equation}\label{SmallEvolution}
\left\lbrace
\begin{array}{rl}
\partial_{s_0} \psi(x,s_0)=& -\nabla\cdot(v\, \psi)(x,s_0)-\Lambda\psi(x,s_0)\\[5mm]
\psi(x,0)=& \psi_0(x)\\[5mm]
div(v)=&0 \quad \mbox{and }\; v\in L^{\infty}([0,T];bmo(\mathbb{R}^n))\quad \mbox{with } \underset{s_0\in [0,T]}{\sup}\; \|v(\cdot,s_0)\|_{bmo}\leq \mu
\end{array}
\right.
\end{equation}
If $\psi_0$ is a small $r$-molecule in the sense of the Definition \ref{DefMolecules} for the local Hardy space $h^\sigma(\mathbb{R}^n)$, then there exists a positive constant $K=K(\mu)$ big enough and a positive constant $\epsilon$ such that for all $0< s_0 \leq\epsilon r$ small we have the following estimates
\begin{eqnarray}
\int_{\mathbb{R}^n}|\psi(x,s_0)||x-x(s_0)|^\omega dx &\leq &(r+Ks_0)^{\omega-\gamma}  \label{SmallConcentration}\\
\|\psi(\cdot, s_0)\|_{L^\infty}&\leq & \frac{1}{\big(r+Ks_0\big)^{n+\gamma}}\label{SmallLinftyevolution}
\end{eqnarray}
This remains true as long as $(r+Ks_0)<1$. The new molecule's center $x(s_0)$ used in formula (\ref{SmallConcentration}) is fixed by 
\begin{equation}\label{Defpointx_0}
\left\lbrace
\begin{array}{rl}
x'(s_0)=& \overline{v}_{B_r}=\frac{1}{|B_r|}\displaystyle{\int_{B_r}}v(y,s_0)dy \qquad \mbox{where }  B_r=B(x(s_0),r)\\[5mm]
x(0)=& x_0.
\end{array}
\right.
\end{equation}  
\end{Theoreme}
\begin{Remarque}
\emph{
\begin{itemize}
\item[1)] The definition of the point $x(s_0)$ given by (\ref{Defpointx_0}) reflects the molecule's center transport using velocity $v$.
\item[2)] Estimates (\ref{SmallConcentration})-(\ref{SmallLinftyevolution}) explain the molecules' deformation following the evolution of the system. Note in particular that if $s_0\longrightarrow 0$ in (\ref{SmallConcentration}) and (\ref{SmallLinftyevolution}) we recover the initial molecular conditions (\ref{Hipo1}) and (\ref{Hipo2}).
\end{itemize}
}
\end{Remarque}
\begin{Corollaire}\label{CorollaireL1}
With inequalities  (\ref{SmallConcentration})-(\ref{SmallLinftyevolution}) of the previous theorem we obtain
\begin{equation}\label{SmallL1evolution}
\|\psi(\cdot, s_0)\|_{L^1}\leq \frac{v_n}{(r+Ks_0)^{\gamma}}
\end{equation} 
where $v_n$ denotes the volume of the $n$-dimensional unit ball. 
\end{Corollaire}
\textit{\textbf{Proof.}}  We write
\begin{eqnarray}
\int_{\mathbb{R}^n}|\psi(x,s_0)|dx&=&\int_{\{|x-x(s_0)|< D\}}|\psi(x,s_0)|dx+\int_{\{|x-x(s_0)|\geq D\}}|\psi(x,s_0)|dx\nonumber\\
&\leq & v_n D^n \|\psi(\cdot, s_0)\|_{L^\infty}+D^{-\omega}\int_{\mathbb{R}}|\psi(x,s_0)||x-x(s_0)|^\omega dx\nonumber
\end{eqnarray}
Now using (\ref{SmallLinftyevolution}) and (\ref{SmallConcentration}) one has:
\begin{eqnarray*}
\int_{\mathbb{R}^n}|\psi(x,s_0)|dx&\leq & v_n \frac{D^n }{(r+Ks_0)^{n+\gamma}}  +D^{-\omega}(r+Ks_0)^{\omega-\gamma}
\end{eqnarray*}
To continue, it is enough to choose correctly the real parameter $D$ to obtain
$$\int_{\mathbb{R}^n}|\psi(x,s_0)|dx\leq v_n\frac{(r+Ks_0)^{(\omega-\gamma)\frac{n}{n+\omega}}}{(r+Ks_0)^{\frac{\omega}{n+\omega}(n+\gamma)}}=\frac{v_n}{(r+Ks_0)^{\gamma}}.$$
\begin{flushright}$\blacksquare$\end{flushright}
\begin{Remarque} \emph{This result shows that the $L^1$ control will be a consequence of the concentration and the height conditions. This corollary also explains the fact that it is enough to treat the case $0<(r+Ks_0)<1$. Indeed, if $(r+Ks_0)=1$, thanks to the bound (\ref{SmallL1evolution}), the $L^1$ control will be trivial then for time $s_0$ and beyond: we only need to apply the maximum principle.}
\end{Remarque}
The proof of Theorem \ref{SmallGeneralisacion} follows the next scheme: we first prove with the Proposition \ref{Propo1} the small concentration condition (\ref{SmallConcentration}); then we will see how this inequality implies the height condition (\ref{SmallLinftyevolution}) which is proved in Proposition \ref{Propo2}.
\begin{Proposition}[Small time Concentration condition]\label{Propo1}
Under the  hypothesis of the Theorem \ref{SmallGeneralisacion}, if $\psi_0$ is a small $r$-molecule, then the solution $\psi(x,s)$ of (\ref{SmallEvolution}) satisfies
\begin{equation*}
\int_{\mathbb{R}^n} |\psi(x,s_0)||x-x(s_0)|^{\omega}dx \leq (r+Ks_0)^{\omega-\gamma}
\end{equation*}
for $x(s_0)\in \mathbb{R}^n$ fixed by the formula (\ref{Defpointx_0}) and with $0\leq s_0\leq \epsilon r$.
\end{Proposition}
\textit{\textbf{Proof}.}
Let us write $\Omega(x-x(s_0))=|x-x(s_0)|^{\omega}$ and $\psi(x)=\psi_+(x)-\psi_-(x)$ where the functions $\psi_{\pm}(x)\geq 0$ have disjoint support. We will note $\psi_\pm(x,s_0)$ solutions of (\ref{SmallEvolution}) with $\psi_\pm(x,0)=\psi_\pm(x)$.\\

At this point, we assume the following positivity principle which is proven in section \ref{Sect_PrincipeMax}: 
\begin{Theoreme}\label{TheoPrincipePos}
Let $1\leq p \leq +\infty$, if initial data $\psi_0\in L^p(\mathbb{R}^n)$ is such that $0\leq \psi_0(x)\leq M$, then the associated solution  $\psi(x,s_0)$ of (\ref{SmallEvolution}) satisfies $0\leq \psi(x,s_0)\leq M$ for all $s_0\in [0,T]$.
\end{Theoreme}
Thus, by linearity and using the above theorem we have that $|\psi(x,s_0)|=|\psi_+(x,s_0)-\psi_-(x,s_0)|\leq \psi_+(x,s_0)+\psi_-(x,s_0)$ and we can write
$$\int_{\mathbb{R}^n}|\psi(x,s_0)|\Omega(x-x(s_0))dx\leq\int_{\mathbb{R}^n}\psi_+(x,s_0)\Omega(x-x(s_0))dx+\int_{\mathbb{R}^n}\psi_-(x,s_0)\Omega(x-x(s_0))dx$$
so we only have to treat one of the integrals on the right side above. We have:
\begin{eqnarray*}
I&=&\left|\partial_{s_0} \int_{\mathbb{R}^n}\Omega(x-x(s_0))\psi_+(x,s_0)dx\right|\\
&=&\left|\int_{\mathbb{R}^n}\partial_{s_0} \Omega(x-x(s_0))\psi_+(x,s_0)+\Omega(x-x(s_0))\left[-\nabla\cdot(v\, \psi_+(x,s_0))-\Lambda\psi_+(x,s_0)\right]dx\right|\\
&=&\left|\int_{\mathbb{R}^n}-\nabla\Omega(x-x(s_0))\cdot x'(s_0)\psi_+(x,s_0)+\Omega(x-x(s_0))\left[-\nabla\cdot(v\, \psi_+(x,s_0))-\Lambda\psi_+(x,s_0)\right]dx\right|
\end{eqnarray*}
Using the fact that $v$ is divergence free, we obtain
\begin{equation*}
I=\left|\int_{\mathbb{R}^n}\nabla\Omega(x-x(s_0))\cdot(v-x'(s_0))\psi_+(x,s_0)-\Omega(x-x(s_0))\Lambda\psi_+(x,s_0)dx\right|.
\end{equation*}
Finally, using the definition of $x'(s_0)$ given in (\ref{Defpointx_0}) and replacing $\Omega(x-x(s_0))$ by $|x-x(s_0)|^{\omega}$ we obtain
\begin{equation}\label{smallEstrella}
I\leq c \underbrace{\int_{\mathbb{R}^n}|x-x(s_0)|^{\omega-1}|v-\overline{v}_{B_r}| |\psi_+(x,s_0)|dx}_{I_1} + c\underbrace{\int_{\mathbb{R}^n}|x-x(s_0)|^{\omega-1}|\psi_+(x,s_0)|dx}_{I_2}.
\end{equation}
We will study separately each of the integrals $I_1$ and $I_2$ in the next lemmas:
\begin{Lemme}\label{Lemme1} For integral $I_1$ above we have the estimate $I_1\leq C \mu \; r^{\omega-1-\gamma}$.
\end{Lemme}
\textit{\textbf{Proof}.} We begin by considering the space $\mathbb{R}^n$ as the union of a ball with dyadic coronas centered on $x(s_0)$, more precisely we set $\mathbb{R}^n=B_r\cup \bigcup_{k\geq 1}E_k$ where
\begin{eqnarray}\label{SmallDecoupage}
B_r&=& \{x\in \mathbb{R}^n: |x-x(s_0)|\leq r\}.\\[5mm]
E_k&=& \{x\in \mathbb{R}^n: r2^{k-1}<|x-x(s_0)|\leq r 2^{k}\} \quad \mbox{for } k\geq 1,\nonumber
\end{eqnarray}
\begin{enumerate}
\item[(i)] \underline{Estimations over the ball $B_r$}. Applying the Hölder inequality on integral $I_1$ we obtain
\begin{eqnarray}
\int_{B_r}|x-x(s_0)|^{\omega-1}|v-\overline{v}_{B_r}| |\psi_+(x,s_0)|dx & \leq &\underbrace{\||x-x(s_0)|^{\omega-1}\|_{L^p(B_r)}}_{(1)} \label{Equa1} \\
& \times &\underbrace{\|v-\overline{v}_{B_r}\|_{L^z(B_r)}}_{(2)}\underbrace{\|\psi_+(\cdot, s_0)\|_{L^q(B_r)}}_{(3)}\nonumber
\end{eqnarray}
where $\frac{1}{p}+\frac{1}{z}+\frac{1}{q}=1$ and $p,z,q> 1$. We treat each of the previous terms separately:
\begin{enumerate}
\item[$\bullet$] First observe that for $1<p<n/(1-\omega)$ we have for the term $(1)$ above:
\begin{equation*}
\||x-x(s_0)|^{\omega-1}\|_{L^p(B_r)}\leq C r^{n/p+\omega-1}.
\end{equation*}
\item[$\bullet$] By hypothesis we have $v(\cdot, s_0)\in bmo$, thus
$$\|v-\overline{v}_{B_r}\|_{L^z(B_r)}\leq C|B_r|^{1/z}\|v(\cdot,s_0)\|_{bmo}.$$ 
since $\underset{s_0\in [0,T]}{\sup}\; \|v(\cdot,s_0)\|_{bmo}\leq \mu$ we write for the term $(2)$
\begin{equation*}
\|v-\overline{v}_{B_r}\|_{L^z(B_r)}\leq C\mu\; r^{n/z}.
\end{equation*}
\item[$\bullet$] Finally for $(3)$ by the maximum principle (\ref{PrincipeMax2}) for $L^q$ norms we have $\|\psi_+(\cdot, s_0)\|_{L^q(B_r)}\leq  \|\psi_+(\cdot, 0)\|_{L^q}$; hence using the fact that $\psi_0$ is an $r$-molecule and remark \ref{Remark2}-2) we obtain
\begin{equation*}
\|\psi_+(\cdot, s_0)\|_{L^q(B_r)}\leq  C \bigg[r^{-\gamma}\bigg]^{1/q}\left[\frac{1}{r^{n+\gamma}}\right]^{1-1/q}.
\end{equation*}
\end{enumerate}
We gather all these inequalities together in order to obtain the following estimation for (\ref{Equa1}):
\begin{equation}\label{Bola1}
\int_{B_r}|x-x(s_0)|^{\omega-1}|v-\overline{v}_{B_r}| |\psi_+(x,s_0)|dx\leq C\mu\; r^{\omega-1-\gamma}.
\end{equation}
\item[(ii)] \underline{Estimations for the dyadic corona $E_k$}. Let us note $I_k$ the integral 
$$I_k=\int_{E_k}|x-x(s_0)|^{\omega-1}|v-\overline{v}_{B_r}| |\psi_+(x,s_0)|dx.$$
Since over $E_k$ we have\footnote{recall that $0<\gamma<\omega<1$.} $|x-x(s_0)|^{\omega-1}\leq C 2^{k(\omega-1)}r^{\omega-1}$ we write
\begin{eqnarray*}
I_k&\leq & C2^{k(\omega-1)}r^{\omega-1}\left(\int_{E_k}|v-\overline{v}_{B_{r2^k}}| |\psi_+(x,s_0)|dx+\int_{E_k}|\overline{v}_{B_r}-\overline{v}_{B_{r2^k}}| |\psi_+(x,s_0)|dx\right)
\end{eqnarray*}
where we noted $B_{r2^k}=B(x(s_0),r2^k)$, then
\begin{eqnarray*}
I \leq C2^{k(\omega-1)}r^{\omega-1}\left(\int_{B_{r2^k}}|v-\overline{v}_{B_{r2^k}}| |\psi_+(x,s_0)|dx+\int_{B_{r2^k}}|\overline{v}_{B_r}-\overline{v}_{B_{r2^k}}| |\psi_+(x,s_0)|dx\right).
\end{eqnarray*}
Now, since $v(\cdot, s_0)\in bmo$, using Proposition \ref{PropoBMO1} we have $|\overline{v}_{B_r}-\overline{v}_{B_{r2^k}}| \leq Ck\|v(\cdot,s_0)\|_{bmo}\leq Ck\mu$ and we write
\begin{eqnarray*}
I_k&\leq &C2^{k(\omega-1)}r^{\omega-1}\left(\int_{B_{r2^k}}|v-\overline{v}_{B_{r2^k}}| |\psi_+(x,s_0)|dx+Ck\mu\|\psi_+(\cdot,s_0)\|_{L^1}\right)\\[5mm]
&\leq & C2^{k(\omega-1)}r^{\omega-1}\left(\|\psi_+(\cdot,s_0)\|_{L^{a_0}} \|v-\overline{v}_{B_{r2^k}}\|_{L^{\frac{a_0}{a_0-1}}} +Ck\mu\;  r^{-\gamma}\right)
\end{eqnarray*}
where we used the Hölder inequality with $1<a_0<\frac{n}{n+(\omega-1)}$ and maximum principle for the last term above. Using again the properties of $bmo$ spaces we have
$$I_k\leq  C2^{k(\omega-1)}r^{\omega-1}\left(\|\psi_+(\cdot,0)\|_{L^1}^{1/a_0}\|\psi_+(\cdot,0)\|_{L^\infty}^{1-1/a_0} |B_{r2^k}|^{1-1/a_0}\|v(\cdot,s)\|_{bmo} +Ck\mu r^{-\gamma}\right).$$
Let us now apply estimates given by hypothesis over $\|\psi_+(\cdot,0)\|_{L^1}$, $\|\psi_+(\cdot,0)\|_{L^\infty}$ and $\|v(\cdot,s_0)\|_{bmo}$  to obtain
$$I_k\leq  C2^{k(n-n/a_0+\omega-1)}r^{\omega-1-\gamma}\mu +C2^{k(\omega-1)}k\mu\;  r^{\omega-1-\gamma}.$$
Since $1<a_0<\frac{n}{n+(\omega-1)}$, we have $n-n/a_0+(\omega-1)<0$, so that, summing over each dyadic corona $E_k$, we have
\begin{equation}\label{Coronak}
\sum_{k\geq 1}I_k\leq C\mu\; r^{\omega-1-\gamma}.
\end{equation}
\end{enumerate}
Finally, gathering together estimations (\ref{Bola1}) and (\ref{Coronak}) we obtain the desired conclusion.\hfill $\blacksquare$\\

\begin{Lemme}\label{Lemme2}
For integral $I_2$ in inequality (\ref{smallEstrella}) we have $I_2\leq C r^{\omega-1-\gamma}$.
\end{Lemme}
\textbf{\textit{Proof.}} As for Lemma \ref{Lemme1}, we consider $\mathbb{R}^n$ as the union of a ball with dyadic coronas centered on $x(s_0)$ (cf. (\ref{SmallDecoupage})). \\

\begin{enumerate}
\item[(i)] \underline{Estimations over the ball $B_r$}. We apply now the Hölder inequality in the integral $I_2$ above with $1<a_1<n/(1-\omega)$ and $\frac{1}{a_1}+\frac{1}{b_1}=1$ in order to obtain
\begin{eqnarray*}
\int_{B_r}|x-x(s_0)|^{\omega-1}|\psi_+(x,s_0)|dx &\leq &\||x-x(s_0)|^{\omega-1}\|_{L^{a_1}(B_r)}\|\psi_+(\cdot, s_0)\|_{L^{b_1}(B_r)}\\
&\leq & Cr^{n/a_1+\omega-1}\|\psi_+(\cdot, 0)\|_{L^1}^{1/b_1}\|\psi_+(\cdot, 0)\|_{L^\infty}^{1-1/b_1}.
\end{eqnarray*}
Using hypothesis over $\|\psi_+(\cdot, 0)\|_{L^1}$ and $\|\psi_+(\cdot, 0)\|_{L^\infty}$ we obtain
\begin{equation}\label{smallBola2}
\int_{B_r}|x-x(s_0)|^{\omega-1}|\psi_+(x,s_0)|dx \leq Cr^{\omega-1-\gamma}.
\end{equation}
\item[(ii)]  \underline{Estimations for the dyadic corona $E_k$}. Here we have
\begin{eqnarray*}
\int_{E_k}|x-x(s_0)|^{\omega-1}|\psi_+(x,s_0)|dx & \leq & C2^{k(\omega-1)}r^{\omega-1}\int_{E_k}|\psi_+(x,s_0)|dx\leq C2^{k(\omega-1)}r^{\omega-1}\|\psi_+(\cdot, s_0)\|_{L^1}\\
&\leq & C2^{k(\omega-1)}r^{\omega-1}\|\psi_+(\cdot, 0)\|_{L^1}\leq C2^{k(\omega-1)}r^{\omega-1-\gamma}
\end{eqnarray*}
It is at this step that the flexibility of molecules is essential. Indeed, in the Definition \ref{DefMolecules} we have fixed $0<\gamma<\omega<1$ so we have $\omega-1<0$ and thus, summing over $k\geq 1$, we obtain
\begin{equation}\label{smallCoronak2}
\sum_{k\geq 1}\int_{E_k}|x-x(s_0)|^{\omega-1}|\psi_+(x,s_0)|dx \leq C r^{\omega-1-\gamma}.
\end{equation}
\end{enumerate}
In order to finish the proof of the Lemma \ref{Lemme2} we glue together estimates (\ref{smallBola2}) and (\ref{smallCoronak2}). \hfill $\blacksquare$\\

Now we continue the proof of the Proposition \ref{Propo1}. Using the Lemmas \ref{Lemme1} and \ref{Lemme2} and getting back to estimate (\ref{smallEstrella}) we have
$$\left|\partial_{s_0} \int_{\mathbb{R}^n}\Omega(x-x(s_0))\psi_+(x,s_0)dx\right| \leq  C(\mu+1)\; r^{\omega-1-\gamma}$$
This last estimation is compatible with the estimate (\ref{SmallConcentration}) for $0\leq s_0\leq \epsilon r$ small enough: just fix $K$ such that
\begin{equation}\label{SmallConstants}
C\left(\mu+1\right)\leq K(\omega-\gamma).
\end{equation}
Indeed, since the time $s_0$ is very small, we can linearize the right-hand side of (\ref{SmallConcentration}) in order to obtain
\begin{equation}\label{FonctionPhi}
\phi=(r+Ks_0)^{\omega-\gamma} \thickapprox r^{\omega-\gamma}\left(1+[K(\omega-\gamma)]\frac{s_0}{r}\right).
\end{equation}
Finally, taking the derivative with respect to $s_0$ in the above expression we have $\phi' \thickapprox r^{\omega-1-\gamma}K(\omega-\gamma)$ and with condition (\ref{SmallConstants}), the Proposition \ref{Propo1} follows. \hfill $\blacksquare$\\

Now we will give a sligthly different proof of the maximum principle of A. C\'ordoba \& D. C\'ordoba. Indeed, the following proof only relies on the concentration condition proved in the lines above. 
\begin{Proposition}[Small time Height condition]\label{Propo2}
Under the  hypothesis of the Theorem \ref{SmallGeneralisacion}, if $\psi(x,s_0)$ satisfies concentration condition (\ref{SmallConcentration}), then we have the next height condition
\begin{equation*}
\|\psi(\cdot, s_0)\|_{L^\infty}  \leq \frac{1}{\left(r+K s_0\right)^{n+\gamma}}.
\end{equation*}
\end{Proposition}
\textit{\textbf{Proof.}} Assume that molecules we are working with are smooth enough. Following an idea of  \cite{Cordoba} (section 4 p.522-523), we will note $\overline{x}$ the point of $\mathbb{R}^n$ such that $\psi(\overline{x},s_0)=\|\psi(\cdot,s_0)\|_{L^\infty}$. Thus we can write
\begin{equation}\label{Infty1}
\frac{d}{ds_0}\|\psi(\cdot,s_0)\|_{L^\infty}\leq -\int_{\mathbb{R}^n}\frac{|\psi(\overline{x},s_0)-\psi(y,s_0)|}{|\overline{x}-y|^{n+1}}dy\leq 0.
\end{equation}
For simplicity, we will assume that $\psi(\overline{x},s_0)$ is positive. Let us consider the corona centered in $\overline{x}$ defined by
$$\mathcal{C}(R_1,R_2)=\{y\in \mathbb{R}^n:R_1\leq|\overline{x}-y|\leq R_2\}$$ 
where $R_2=\rho R_1$ with $\rho >2$ and where $R_1$ will be fixed later. Then:
\begin{equation*}
\int_{\mathbb{R}^n}\frac{\psi(\overline{x},s_0)-\psi(y,s_0)}{|\overline{x}-y|^{n+1}}dy\geq \int_{\mathcal{C}(R_1,R_2)}\frac{\psi(\overline{x},s_0)-\psi(y,s_0)}{|\overline{x}-y|^{n+1}}dy.
\end{equation*}
Define the sets $B_1$ and $B_2$ by $B_1=\{y\in \mathcal{C}(R_1,R_2): \psi(\overline{x},s_0)-\psi(y,s_0)\geq \frac{1}{2}\psi(\overline{x},s_0)\}$ and $B_2=\{y\in \mathcal{C}(R_1,R_2): \psi(\overline{x},s_0)-\psi(y,s_0)< \frac{1}{2}\psi(\overline{x},s_0)\}$ such that $\mathcal{C}(R_1,R_2)=B_1\cup B_2$. We obtain the inequalities
\begin{eqnarray*}
\int_{\mathcal{C}(R_1,R_2)}\frac{\psi(\overline{x},s_0)-\psi(y,s_0)}{|\overline{x}-y|^{n+1}}dy \geq \int_{B_1}\frac{\psi(\overline{x},s_0)-\psi(y,s_0)}{|\overline{x}-y|^{n+1}}dy \geq \frac{\psi(\overline{x},s_0)}{2R_2^{n+1}}|B_1|=\frac{\psi(\overline{x},s_0)}{2R_2^{n+1}}\left(|\mathcal{C}(R_1,R_2)|-|B_2|\right).
\end{eqnarray*}
Since $R_2=\rho R_1$ one has
\begin{equation}\label{Infty3}
\int_{\mathcal{C}(R_1,R_2)}\frac{\psi(\overline{x},s_0)-\psi(y,s_0)}{|\overline{x}-y|^{n+1}}dy\geq \frac{\psi(\overline{x},s_0)}{2\rho^{n+1}R_1^{n+1}}\bigg(v_n(\rho^n -1)R_1^n-|B_2|\bigg)
\end{equation}
where $v_n$ denotes the volume of the $n$-dimensional unit ball. Now, we will estimate the quantity $|B_2|$ in terms of $\psi(\overline{x},s_0)$ and $R_1$ with the next lemma.
\begin{Lemme}\label{LemmaB2}
For the set $B_2$ we have the following estimations
\begin{enumerate}
\item[1)] if $|\overline{x}-x(s_0)|>2R_2$ then $(r+Ks_0)^{\omega-\gamma}C_1\psi(\overline{x},s_0)^{-1}R_1^{-\omega}\geq |B_2|$.\\

\item[2)] if $|\overline{x}-x(s_0)|<R_1/2$ then $(r+Ks_0)^{\omega-\gamma}C_1\psi(\overline{x},s_0)^{-1}R_1^{-\omega}\geq |B_2|$.\\

\item[3)] if $R_1/2\leq |\overline{x}-x(s_0)|\leq 2R_2$ then $(r+Ks_0)^{\frac{\omega}{2}-\frac{\gamma}{2}}(C_2 R_1^{n-\omega}\psi(\overline{x},s_0)^{-1})^{1/2}\geq |B_2|$.
\end{enumerate}
\end{Lemme}
Recall that for the molecule's center $x_0\in \mathbb{R}^n$ we noted its transport by $x(s_0)$ which is defined by formula (\ref{Defpointx_0}).\\[5mm] \textit{\textbf{Proof.}} For all these estimates, our starting point is the concentration condition (\ref{SmallConcentration}):
\begin{eqnarray}
(r+Ks_0)^{\omega-\gamma}\geq \int_{\mathbb{R}^n}|\psi(y,s_0)||y-x(s_0)|^{\omega}dy \geq \int_{B_2}|\psi(y,s_0)||y-x(s_0)|^{\omega}dy \geq \frac{\psi(\overline{x},s_0)}{2}\int_{B_2}|y-x(s_0)|^{\omega}dy.\label{EstimationB2}
\end{eqnarray}
We just need to estimate the last integral following the cases given by the lemma. The first two cases are very similar. Indeed, if $|\overline{x}-x(s_0)|>2R_2$  then we have
$$\underset{y\in B_2\subset \mathcal{C}(R_1,R_2)}{\min}|y-x(s_0)|^{\omega}\geq R_2^{\omega}=\rho^{\omega}R_1^{\omega}$$
while for the second case, if $|\overline{x}-x(s_0)|<R_1/2$, one has
$$\underset{y\in B_2\subset \mathcal{C}(R_1,R_2)}{\min}|y-x(s_0)|^{\omega}\geq \frac{R_1^{\omega}}{2^\omega}.$$
Applying these results to (\ref{EstimationB2}) we obtain $(r+Ks_0)^{\omega-\gamma}\geq \frac{\psi(\overline{x},s_0)}{2} \rho^{\omega} R_1^{\omega}|B_2|$ and  $(r+Ks_0)^{\omega-\gamma}\geq \frac{\psi(\overline{x},s_0)}{2} \frac{R_1^{\omega}}{2^\omega}|B_2|$, and since $\rho>2$ we have the desired estimate
\begin{equation*}
\frac{(r+Ks_0)^{\omega-\gamma}C_1}{\psi(\overline{x},s_0) R_1^{\omega}} \geq \frac{2(r+Ks_0)^{\omega-\gamma}}{\rho^{\omega}\psi(\overline{x},s_0) R_1^{\omega}} \geq |B_2|
\end{equation*}
with $C_1=2^{1+\omega}$. For the last case, since $R_1/2\leq |\overline{x}-x(s_0)|\leq 2R_2$ we can write using the Cauchy-Schwarz inequality
\begin{equation}\label{HolderInver}
\int_{B_2}|y-x(s_0)|^{\omega}dy\geq |B_2|^2\left(\int_{B_2}|y-x(s_0)|^{-\omega}dy\right)^{-1}
\end{equation}
Now, observe that in this case we have $B_2\subset B(x(s_0), 5R_2)$ and then
$$\int_{B_2}|y-x(s_0)|^{-\omega}dy\leq \int_{B(x(s_0), 5 R_2)}|y-x(s_0)|^{-\omega}dy\leq v_n (5\rho R_1)^{n-\omega}.$$
Getting back to (\ref{HolderInver}) we obtain
$$\int_{B_2}|y-x(s_0)|^{\omega}dy\geq |B_2|^2 v_n^{-1} (5 \rho R_1)^{-n+\omega}$$
We use this estimate in (\ref{EstimationB2}) to obtain
\begin{equation*}
(r+Ks_0)^{\frac{\omega}{2}-\frac{\gamma}{2}}\frac{C_2 R_1^{n/2-\omega/2}}{\psi(\overline{x},s_0)^{1/2}}\geq |B_2|,
\end{equation*}
where $C_2=(2\times 5^{n-\omega} v_n\rho^{n-\omega})^{1/2}$. The lemma is proven.\hfill $\blacksquare$\\

With this lemma at our disposal we can write
\begin{enumerate}
\item[(i)] if $|\overline{x}-x(s_0)|>2R_2$ or $|\overline{x}-x(s_0)|<R_1/2$ then
\begin{equation*}
\int_{\mathcal{C}(R_1,R_2)}\frac{\psi(\overline{x},s_0)-\psi(y,s_0)}{|\overline{x}-y|^{n+1}}dy\geq  \frac{\psi(\overline{x},s_0)}{2\rho^{n+1}R_1^{n+1}}\bigg(v_n(\rho^n -1)R_1^n-\frac{C_1(r+Ks_0)^{\omega-\gamma}}{\psi(\overline{x},s_0)} R_1^{-\omega}\bigg)
\end{equation*}
\item[(ii)] if $R_1/2\leq |\overline{x}-x(s_0)|\leq 2R_2$
\begin{equation*}
\int_{\mathcal{C}(R_1,R_2)}\frac{\psi(\overline{x},s_0)-\psi(y,s_0)}{|\overline{x}-y|^{n+1}}dy\geq  \frac{\psi(\overline{x},s_0)}{2\rho^{n+1}R_1^{n+1}}\bigg(v_n(\rho^n -1)R_1^n-\frac{C_2 (r+Ks_0)^{\frac{\omega}{2}-\frac{\gamma}{2}}R_1^{n/2-\omega/2}}{\psi(\overline{x},s_0)^{1/2}}\bigg)
\end{equation*}
\end{enumerate}
Now, if we set $R_1=(r+Ks_0)^{\frac{(\omega-\gamma)}{n+\omega}}\psi(\overline{x},s_0)^{\frac{-1}{n+\omega}}$ and if $\rho$ is big enough, we obtain for cases (i) and (ii) the following estimate for (\ref{Infty3}):
\begin{equation*}
\int_{\mathcal{C}(R_1,R_2)}\frac{\psi(\overline{x},s_0)-\psi(y,s_0)}{|\overline{x}-y|^{n+1}}dy\geq C (r+Ks_0)^{-\frac{(\omega-\gamma)}{n+\omega}} \psi(\overline{x},s_0)^{1+\frac{1}{n+\omega}}
\end{equation*}
where $C=C(n,\rho)=\frac{v_n (\rho^n-1)-\sqrt{2v_n}(5\rho)^{\frac{n-\omega}{2}}}{2\rho^{n+1}}<1$ is a small positive constant. Hence, and for all possible cases considered before, we have the next estimate for (\ref{Infty1}):
$$\frac{d}{ds_0}\|\psi(\cdot,s_0)\|_{L^\infty}\leq-C(r+Ks_0)^{-\frac{(\omega-\gamma)}{n+\omega}}\|\psi(\cdot,s_0)\|_{L^\infty}^{1+\frac{1}{n+\omega}}.$$
In order to solve this problem it is enough to remark that if $\|\psi(\cdot,s_0)\|_{L^\infty}\leq (r+Ks_0)^{-(n+\gamma)}$ then $\|\psi(\cdot,s_0)\|_{L^\infty}$ satisfies the previous estimate. Indeed, recalling that $(r+Ks_0)<1$, we have
\begin{eqnarray*}
\frac{d}{ds_0}\|\psi(\cdot,s_0)\|_{L^\infty}&\leq & -K(n+\gamma)(r+Ks_0)^{-(n+\gamma)-1}\\
&\leq &-C(r+Ks_0)^{-\frac{(\omega-\gamma)}{n+\omega}}(r+Ks_0)^{-(n+\gamma)(1+\frac{1}{n+\omega})}\\
&\leq &-C(r+Ks_0)^{-\frac{(\omega-\gamma)}{n+\omega}}\|\psi(\cdot,s_0)\|_{L^\infty}^{1+\frac{1}{n+\omega}}.
\end{eqnarray*}
Furthermore, this solution is unique. \hfill $\blacksquare$
\subsection{Molecule's evolution: Second step}\label{SecEvolMol2}
In the previous section we have obtained deformed molecules after a very small time $s_0$. The next theorem shows us how to obtain similar profiles in the inputs and the outputs in order to perform an iteration in time. 
\begin{Theoreme}\label{Generalisacion} Set $\gamma$ and $\omega$ two real numbers such that $0<\gamma<\omega<1$. Let $0< s_1\leq T$ and let $\psi(x,s_1)$ be a solution of the problem
\begin{equation}\label{Evolution}
\left\lbrace
\begin{array}{rl}
\partial_{s_1} \psi(x,s_1)=& -\nabla\cdot(v\, \psi)(x,s_1)-\Lambda\psi(x,s_1)\\[5mm]
\psi(x,0)=& \psi(x,s_0)  \qquad \qquad \mbox{with } s_0>0\\[5mm]
div(v)=&0 \quad \mbox{and }\; v\in L^{\infty}([0,T];bmo(\mathbb{R}^n))\quad \mbox{with } \underset{s_1\in [s_0,T]}{\sup}\; \|v(\cdot,s_1)\|_{bmo}\leq \mu
\end{array}
\right.
\end{equation}
If $\psi(x,s_0)$ satisfies the three following conditions
\begin{eqnarray}
\int_{\mathbb{R}^n}|\psi(x,s_0)||x-x(s_0)|^\omega dx &\leq &(r+Ks_0)^{\omega-\gamma}\nonumber\\
\|\psi(\cdot, s_0)\|_{L^\infty}&\leq & \frac{1}{\left(r+ Ks_0\right)^{n+\gamma}}\nonumber\\
\|\psi(\cdot, s_0)\|_{L^1} &\leq & \frac{v_n}{(r+Ks_0)^{\gamma}}\label{Hipoth}
\end{eqnarray}
where $K=K(\mu)$ is given by (\ref{SmallConstants}) and $s_0$ is such that $0<(r+Ks_0)<1$. Then for all $0< s_1\leq\epsilon r$ small, we have the following estimates
\begin{eqnarray}
\int_{\mathbb{R}^n}|\psi(x,s_1)||x-x(s_1)|^\omega dx &\leq &(r +K(s_0+s_1))^{\omega-\gamma}  \label{Concentration2}\\
\|\psi(\cdot,s_1)\|_{L^\infty}&\leq & \frac{1}{\left(r+K(s_0+s_1)\right)^{n+\gamma}}\label{Linftyevolutionnotsmalltime}\\
\|\psi(\cdot,s_1)\|_{L^1} &\leq & \frac{v_n}{\big(r+K(s_0+s_1)\big)^{\gamma}} \label{L1evolutionsmalltime}
\end{eqnarray}
\end{Theoreme}
\begin{Remarque}\label{RemarkStep1}
\emph{
\begin{itemize}
\item[1)] The $L^1$ bound in (\ref{Hipoth}) is given as hypothesis in order to simplify the exposition: it is not a part of the definition of the molecules. 
\item[2)] As for the Theorem \ref{SmallGeneralisacion}, we only need to study the case when $0<(r+K(s_0+s_1))<1$. Indeed, if $(r+K(s_0+s_1))=1$ and since $s_1$ is small, we obtain immediately that $\|\psi(\cdot, s_1)\|_{L^1}<+\infty$.
\item[3)] The new molecule's center $x(s_1)$ used in formula (\ref{Concentration2}) is fixed by 
\begin{equation}\label{Defpointx_s}
\left\lbrace
\begin{array}{rl}
x'(s_1)=& \overline{v}_{B_{f_1}}=\frac{1}{|B_{f_1}|}\displaystyle{\int_{B_{f_1}}}v(y,s_1)dy\\[5mm]
x(0)=& x(s_0).
\end{array}
\right.
\end{equation}  
And here we noted $B_{f_1}=B(x(s_1),f_1)$ with $f$ a real valued function given by 
\begin{equation}\label{DefiFunctionF}
f_1= \big(r+Ks_0\big).
\end{equation}
\end{itemize}
}
\end{Remarque}
To prove this theorem we will follow the same scheme as before: first we prove the concentration condition (\ref{Concentration2}) in the Proposition \ref{Propo4}. With this estimate at hand we will control the $L^\infty$ decay in Proposition \ref{Propo5} and then we will obtain the suitable $L^1$ control in Proposition \ref{Propo6}.
\begin{Proposition}[Concentration condition]\label{Propo4}
Under the  hypothesis of the Theorem \ref{Generalisacion}, if $\psi(\cdot,s_0)$ is an initial data then the solution $\psi(x,s_1)$ of (\ref{Evolution}) satisfies
\begin{equation*}
\int_{\mathbb{R}^n} |\psi(x,s_1)||x-x(s_1)|^{\omega}dx \leq (r+K(s_0+s_1))^{\omega-\gamma}
\end{equation*}
for $x(s_1)\in \mathbb{R}^n$ fixed by formula (\ref{Defpointx_s}), with $0\leq s_1\leq \epsilon r$. 
\end{Proposition}
\textit{\textbf{Proof}.} The calculations are very similar of those of the Proposition \ref{Propo1}: the only diference stems from the initial data and the definition of the center $x(s_1)$. So, let us write $\Omega(x-x(s_1))=|x-x(s_1)|^{\omega}$ and $\psi(x)=\psi_+(x)-\psi_-(x)$ where the functions $\psi_{\pm}(x)\geq 0$ have disjoint support. 
Thus, by linearity and using the positivity theorem we have  
$$|\psi(x,s_1)|=|\psi_+(x,s_1)-\psi_-(x,s_1)|\leq \psi_+(x,s_1)+\psi_-(x,s_1)$$ and we can write
$$\int_{\mathbb{R}^n}|\psi(x,s_1)|\Omega(x-x(s_1))dx\leq\int_{\mathbb{R}^n}\psi_+(x,s_1)\Omega(x-x(s_1))dx+\int_{\mathbb{R}^n}\psi_-(x,s_1)\Omega(x-x(s_1))dx$$
so we only have to treat one of the integrals on the right-hand side above. We have:
\begin{eqnarray*}
I&=&\left|\partial_{s_1} \int_{\mathbb{R}^n}\Omega(x-x(s_1))\psi_+(x,s_1)dx\right|\\
&=&\left|\int_{\mathbb{R}^n}\partial_{s_1} \Omega(x-x(s_1))\psi_+(x,s_1)+\Omega(x-x(s_1))\left[-\nabla\cdot(v\, \psi_+(x,s_1))-\Lambda\psi_+(x,s_1)\right]dx\right|\\
&=&\left|\int_{\mathbb{R}^n}-\nabla\Omega(x-x(s_1))\cdot x'(s_1)\psi_+(x,s_1)+\Omega(x-x(s_1))\left[-\nabla\cdot(v\, \psi_+(x,s_1))-\Lambda\psi_+(x,s_1)\right]dx\right|
\end{eqnarray*}
Using the fact that $v$ is divergence free, we obtain
\begin{equation*}
I=\left|\int_{\mathbb{R}^n}\nabla\Omega(x-x(s_1))\cdot(v-x'(s_1))\psi_+(x,s_1)-\Omega(x-x(s_1))\Lambda\psi_+(x,s_1)dx\right|.
\end{equation*}
Finally, using the definition of $x'(s_1)$ given in (\ref{Defpointx_s}) and replacing $\Omega(x-x(s_1))$ by $|x-x(s_1)|^{\omega}$ we obtain
\begin{equation}\label{Estrella}
I\leq c \underbrace{\int_{\mathbb{R}^n}|x-x(s_1)|^{\omega-1}|v-\overline{v}_{B_{f_1}}| |\psi_+(x,s_1)|dx}_{I_1} + c\underbrace{\int_{\mathbb{R}^n}|x-x(s_1)|^{\omega-1}|\psi_+(x,s_1)|dx}_{I_2}.
\end{equation}
Again, we will study separately each of the integrals $I_1$ and $I_2$ in the next lemmas:
\begin{Lemme}\label{Lemme41} For integral $I_1$ we have the estimate $I_1\leq C \mu\big(r+Ks_0\big)^{\omega-\gamma-1}$.
\end{Lemme}
\textit{\textbf{Proof}.} We begin by considering the space $\mathbb{R}^n$ as the union of a ball with dyadic coronas centered on $x(s_1)$, more precisely we set $\mathbb{R}^n=B_{f_1}\cup\bigcup_{k\geq 1}E_k$ where
\begin{eqnarray}\label{Decoupage}
B_{f_1}&=& \{x\in \mathbb{R}^n: |x-x(s_1)|\leq f_1\},\\[5mm]
E_k&=& \{x\in \mathbb{R}^n: f_12^{k-1}<|x-x(s_1)|\leq f_1 2^{k}\} \quad \mbox{for } k\geq 1.\nonumber
\end{eqnarray}
\begin{enumerate}
\item[(i)] \underline{Estimations over the ball $B_{f_1}$}. Applying Hölder inequality on integral $I_1$ we obtain
\begin{eqnarray*}
I_{1,B_{f_1}}=\int_{B_{f_1}}|x-x(s_1)|^{\omega-1}|v-\overline{v}_{B_{f_1}}| |\psi_+(x,s_1)|dx & \leq &\underbrace{\||x-x(s_1)|^{\omega-1}\|_{L^p(B_{f_1})}}_{(1)} \\
& \times &\underbrace{\|v-\overline{v}_{B_{f_1}}\|_{L^z(B_{f_1})}}_{(2)}\underbrace{\|\psi_+(\cdot, s_1)\|_{L^q(B_{f_1})}}_{(3)}
\end{eqnarray*}
where $\frac{1}{p}+\frac{1}{z}+\frac{1}{q}=1$ and $p,z,q> 1$. 
\begin{enumerate}
\item[$\bullet$] Observe that for $1<p<n/(1-\omega)$ we have
\begin{equation*}
\||x-x(s_1)|^{\omega-1}\|_{L^p(B_{f_1})}\leq C f_1^{n/p+\omega-1}.
\end{equation*}
\item[$\bullet$] We have $v(\cdot, s_1)\in bmo$, thus $\|v-\overline{v}_{B_{f_1}}\|_{L^z(B_{f_1})}\leq C|B_{f_1}|^{1/z}\|v(\cdot,s_1)\|_{bmo}$, since $\underset{s_1\in [s_0,T]}{\sup}\; \|v(\cdot,s_1)\|_{bmo}\leq \mu$ we write
\begin{equation*}
\|v-\overline{v}_{B_{f_1}}\|_{L^z(B_{f_1})}\leq C f_1^{n/z}\mu.
\end{equation*}
\item[$\bullet$] Finally, by the maximum principle for $L^q$ norms we have:
\begin{equation*}
\|\psi_+(\cdot, s_1)\|_{L^q(B_{f_1})}\leq  \|\psi(\cdot, s_0)\|_{L^1}^{1/q}\|\psi(\cdot, s_0)\|_{L^\infty}^{1-1/q}.
\end{equation*}
\end{enumerate}
We gather all these inequalities in order to obtain the following estimation for $I_{1,B_{f_1}}$:
\begin{equation*}
I_{1,B_{f_1}}=\int_{B_{f_1}}|x-x(s_1)|^{\omega-1}|v-\overline{v}_{B_{f_1}}| |\psi_+(x,s_1)|dx\leq C\mu f_1^{n(1-1/q)+\omega-1}\|\psi(\cdot, s_0)\|_{L^1}^{1/q}\|\psi(\cdot, s_0)\|_{L^\infty}^{1-1/q}.
\end{equation*}
\item[(ii)] \underline{Estimations for the dyadic corona $E_k$}. Let us note $I_{1,E_k}$ the integral
$$I_{1,E_k}=\int_{E_k}|x-x(s_1)|^{\omega-1}|v-\overline{v}_{B_{f_1}}| |\psi_+(x,s_1)|dx.$$
Since over $E_k$ we have $|x-x(s_1)|^{\omega-1}\leq C 2^{k(\omega-1)}f_1^{\omega-1}$ we write
\begin{eqnarray*}
I_{1,E_k}&\leq & C2^{k(\omega-1)}f^{\omega-1}\left(\int_{E_k}|v-\overline{v}_{B(f_12^k)}| |\psi_+(x,s_1)|dx+\int_{E_k}|\overline{v}_{B_{f_1}}-\overline{v}_{B(f_12^k)}| |\psi_+(x,s_1)|dx\right)\\
&\leq & C2^{k(\omega-1)}f_1^{\omega-1}\left(\int_{B(f_12^k)}|v-\overline{v}_{B(f_12^k)}| |\psi_+(x,s_1)|dx+\int_{B(f_12^k)}|\overline{v}_{B_{f_1}}-\overline{v}_{B(f_12^k)}| |\psi_+(x,s_1)|dx\right),
\end{eqnarray*}
where $B(f_12^k)=B(x(s_1), f_1 2^k)$. Now, since $v(\cdot, s_1)\in bmo$, using Proposition \ref{PropoBMO1} we have $|\overline{v}_{B_{f_1}}-\overline{v}_{B(f_12^k)}| \leq Ck\|v(\cdot,s_1)\|_{bmo}\leq Ck\mu$ and we can write
\begin{eqnarray*}
I_{1,E_k}&\leq &C2^{k(\omega-1)}f_1^{\omega-1}\left(\int_{B(f_12^k)}|v-\overline{v}_{B(f_12^k)}| |\psi_+(x,s_1)|dx+Ck\mu\|\psi_+(\cdot,s_1)\|_{L^1}\right)\\[5mm]
&\leq & C2^{k(\omega-1)}f_1^{\omega-1}\left(\|\psi_+(\cdot,s_1)\|_{L^{a_0}} \|v-\overline{v}_{B(f_12^k)}\|_{L^{\frac{a_0}{a_0-1}}} +Ck\mu\;\|\psi_+(\cdot,s_0)\|_{L^1}\right)
\end{eqnarray*}
where we used Hölder inequality with $1<a_0<\frac{n}{n+(\omega-1)}$ and maximum principle for the last term above. Using again the properties of $bmo$ spaces we have
$$I_{1,E_k}\leq  C2^{k(\omega-1)}f_1^{\omega-1}\left(\|\psi_+(\cdot,s_0)\|_{L^1}^{1/a_0}\|\psi_+(\cdot,s_0)\|_{L^\infty}^{1-1/a_0} |B(f_12^k)|^{1-1/a_0}\|v(\cdot,s_1)\|_{bmo} +Ck\mu\|\psi(\cdot,s_0)\|_{L^1}\right).$$
Since $\|v(\cdot,s_1)\|_{bmo}\leq \mu$ and since $1<a_0<\frac{n}{n+(\omega-1)}$, we have $n(1-1/a_0)+(\omega-1)<0$, so that, summing over each dyadic corona $E_k$, we obtain
\begin{equation*}
\sum_{k\geq 1}I_{1,E_k}\leq C \mu\left(f_1^{n(1-1/a_0)+\omega-1}\|\psi(\cdot,s_0)\|_{L^1}^{1/a_0}\|\psi(\cdot,s_0)\|_{L^\infty}^{1-1/a_0}+f_1^{\omega-1}\|\psi(\cdot,s_0)\|_{L^1}\right).
\end{equation*}
\end{enumerate}
We finally obtain the following inequalities:
\begin{eqnarray}
I_1&=& I_{1,B_{f_1}}+\sum_{k\geq 1}I_{1,E_k}\label{FinalEstimateI_1}\\
&\leq &C\mu \underbrace{f_1^{n(1-1/q)+\omega-1}\|\psi(\cdot, s_0)\|_{L^1}^{1/q}\|\psi(\cdot, s_0)\|_{L^\infty}^{1-1/q}}_{(a)}\nonumber\\
& &+ C \mu \left(\underbrace{f_1^{n(1-1/a_0)+\omega-1}\|\psi(\cdot,s_0)\|_{L^1}^{1/a_0}\|\psi(\cdot,s_0)\|_{L^\infty}^{1-1/a_0}}_{(b)} +\underbrace{f_1^{\omega-1}\|\psi(\cdot,s_0)\|_{L^1}}_{(c)}\right)\nonumber
\end{eqnarray}
Now we will prove that each of the terms $(a)$, $(b)$ and $(c)$ above is bounded by the quantity $\big(r+Ks_0\big)^{\omega-\gamma-1}$. Indeed:
\begin{itemize}
\item for the first term (a) by the hypothesis on the initial data $\psi(\cdot,s_0)$ and the definition of $f_1$ given in (\ref{DefiFunctionF}) we have:
\begin{eqnarray*}
f_1^{n(1-1/q)+\omega-1}\|\psi(\cdot, s_0)\|_{L^1}^{1/q}\|\psi(\cdot, s_0)\|_{L^\infty}^{1-1/q}& \leq & \big(r+Ks_0\big)^{n(1-1/q)+\omega-1}\\
& & \times \left[(r+Ks_0)^{-\gamma}\right]^{1/q}\left[(r+Ks_0)^{-(n+\gamma)}\right]^{1-1/q}\\
& \leq & \big(r+Ks_0\big)^{\omega-\gamma-1}.
\end{eqnarray*}
\item For the second term (b) we have, by the same arguments:
\begin{eqnarray*}
f_1^{n(1-1/a_0)+\omega-1}\|\psi(\cdot,s_0)\|_{L^1}^{1/a_0}\|\psi(\cdot,s_0)\|_{L^\infty}^{1-1/a_0} &\leq & \big(r+Ks_0\big)^{n(1-1/a_0)+\omega-1}\\
& & \times \left[(r+Ks_0)^{-\gamma}\right]^{1/a_0}\left[(r+Ks_0)^{-(n+\gamma)}\right]^{1-1/a_0}\\
& \leq & \big(r+Ks_0\big)^{\omega-\gamma-1}.
\end{eqnarray*}
\item Finally, for the last term (c) we write
\begin{eqnarray*}
f_1^{\omega-1}\|\psi(\cdot,s_0)\|_{L^1}&\leq &f_1^{\omega-1}(r+Ks_0)^{-\gamma}= (r+Ks_0)^{\omega-\gamma-1}\\
\end{eqnarray*}
\end{itemize}
Gathering these estimates on $(a), (b)$ and $(c)$, and getting back to (\ref{FinalEstimateI_1}) we finally obtain $I_1\leq C\mu \big(r+Ks_0\big)^{\omega-\gamma-1}$. \hfill$\blacksquare$\\
\begin{Lemme}\label{Lemme42}
For integral $I_2$ in inequality (\ref{Estrella}) we have $I_2\leq C\big(r+Ks_0\big)^{\omega-\gamma-1}$.
\end{Lemme}
\textbf{\textit{Proof.}} As for the Lemma \ref{Lemme41}, we consider $\mathbb{R}^n$ as the union of a ball with dyadic coronas centered on $x(s_1)$ (cf. (\ref{Decoupage})). 
\begin{enumerate}
\item[(i)] \underline{Estimations over the ball $B_{f_1}$}. Applying Hölder inequality with $1<a_1<n/(1-\omega)$ and maximum principle we have
\begin{eqnarray*}
I_{2,B_{f_1}}=\int_{B_{f_1}}|x-x(s_1)|^{\omega-1}|\psi_+(x,s_1)|dx &\leq &\||x-x(s_1)|^{\omega-1}\|_{L^{a_1}(B_{f_1})}\|\psi_+(\cdot, s_1)\|_{L^{b_1}(B_{f_1})}\\
&\leq & Cf_1^{n/a_1+\omega-1}\|\psi_+(\cdot, s_0)\|_{L^1}^{1-1/a_1}\|\psi_+(\cdot, s_0)\|_{L^\infty}^{1/a_1}.
\end{eqnarray*}
Thus, we can write:
\begin{equation}\label{Bola2}
I_{2,B_{f_1}}\leq Cf_1^{n/a_1+\omega-1}\|\psi_+(\cdot, s_0)\|_{L^1}^{1-1/a_1}\|\psi_+(\cdot, s_0)\|_{L^\infty}^{1/a_1}.
\end{equation}
\item[(ii)]  \underline{Estimations for the dyadic corona $E_k$}. Here we have
\begin{eqnarray*}
I_{2,E_k}=\int_{E_k}|x-x(s_1)|^{\omega-1}|\psi_+(x,s_1)|dx & \leq & C2^{k(\omega-1)}f_1^{\omega-1}\int_{E_k}|\psi_+(x,s_1)|dx\leq C2^{k(\omega-1)}f_1^{\omega-1}\|\psi_+(\cdot, s_1)\|_{L^1}\\
&\leq & C2^{k(\omega-1)}f_1^{\omega-1}\|\psi(\cdot, s_0)\|_{L^1}
\end{eqnarray*}
Since $0<\gamma<\omega<1$ we have $\omega-1<0$ and thus, summing over $k\geq 1$, we obtain
\begin{equation*}
\sum_{k\geq 1}I_{2,E_k}=\sum_{k\geq 1}\int_{E_k}|x-x(s_1)|^{\omega-1}|\psi_+(x,s_1)|dx \leq Cf_1^{\omega-1}\|\psi(\cdot, s_0)\|_{L^1}.
\end{equation*}
We finally get
\begin{equation}\label{Coronak2}
\sum_{k\geq 1}I_{2,E_k}\leq C f_1^{\omega-1}\|\psi(\cdot, s_0)\|_{L^1}.
\end{equation}
\end{enumerate}
To finish the proof of the Lemma \ref{Lemme42}  we glue together (\ref{Bola2}) and (\ref{Coronak2}) and we obtain
$$I_2=I_{2,B_{f_1}}+\sum_{k\geq 1}I_{2,E_k}\leq C\left(\underbrace{f_1^{n/a_1+\omega-1}\|\psi_+(\cdot, s_0)\|_{L^1}^{1-1/a_1}\|\psi_+(\cdot, s_0)\|_{L^\infty}^{1/a_1}}_{(d)}+\underbrace{f_1^{\omega-1}\|\psi(\cdot, s_0)\|_{L^1}}_{(e)}\right)$$
Now, we prove that the quantities $(d)$ and $(e)$ can be bounded by $\big(r+Ks_0\big)^{\omega-\gamma-1}$.
\begin{itemize}
\item For the term $(d)$ we write
\begin{eqnarray*}
f_1^{n/a_1+\omega-1}\|\psi(\cdot, s_0)\|_{L^1}^{1-1/a_1}\|\psi(\cdot, s_0)\|_{L^\infty}^{1/a_1}& \leq & (r+Ks_0)^{n/a_1+\omega-1}\\
& & \times \left[(r+Ks_0)^{-\gamma}\right]^{1-1/a_1}\left[(r+Ks_0)^{-(n+\gamma)}\right]^{1/a_1}\\
& \leq & \big(r+Ks_0\big)^{\omega-\gamma-1}
\end{eqnarray*}
\item To treat the term $(e)$ it is enough to apply the same arguments used to prove the part $(c)$ above. 
\end{itemize}
Finally, we obtain $I_2=I_{2,B_{f_1}}+\displaystyle{\sum_{k\geq 1}}I_{2, E_k}\leq  C\big(r+Ks_0\big)^{\omega-\gamma-1}$ and the Lemma \ref{Lemme42} is proven.\hfill$\blacksquare$\\

Now we continue the proof of the Proposition \ref{Propo4}. Using Lemmas \ref{Lemme41} and \ref{Lemme42} and getting back to the estimate (\ref{Estrella}) we have
\begin{equation}\label{FinalEstimate}
\left|\partial_{s_1} \int_{\mathbb{R}^n}\Omega(x-x(s_1))\psi_+(x,s_1)dx\right| \leq  C \left(\mu+1\right)\big(r+Ks_0\big)^{\omega-\gamma-1}
\end{equation}
This estimation is compatible with the estimate (\ref{Concentration2}) for $0\leq s_1\leq \epsilon r$ small enough. Indeed, we can write
\begin{eqnarray*}
\phi&=&\left(r+K(s_0+s_1)\right)^{\omega-\gamma}
\end{eqnarray*}
and we linearize this expression with respect to $s_1$:
$$\phi \thickapprox (r+Ks_0)^{\omega-\gamma}\left(1+K(\omega-\gamma)\frac{s_1}{(r+Ks_0)}\right).$$
Taking the derivative of $\phi$ with respect to $s_1$ we have $\phi' \thickapprox K(\omega-\gamma) \big(r+Ks_0\big)^{\omega-\gamma-1}$ and with the condition (\ref{SmallConstants}) on $K(\omega-\gamma)$ we obtain that (\ref{FinalEstimate}) is bounded by $\phi'$ and the Proposition \ref{Propo4} follows. \hfill $\blacksquare$\\

Now we write down the maximum principle for a small time $s_1$ but with a initial condition $\psi(\cdot,s_0)$, with $s_0>0$.
\begin{Proposition}[Height condition]\label{Propo5}
Under the  hypothesis of the Theorem \ref{Generalisacion}, if $\psi(x,s_1)$ satisfies concentration condition (\ref{Concentration2}), then we have the next height condition
\begin{equation*}
\|\psi(\cdot, s_1)\|_{L^\infty}  \leq \frac{1}{\left(r+K(s_0+s_1)\right)^{n+\gamma}}.
\end{equation*}
\end{Proposition}
\textit{\textbf{Proof.}} The proof follows essentially the same lines of the Proposition \ref{Propo2}. Indeed, since we have that concentration condition (\ref{Concentration2}) is bounded by $\big(r+K(s_0+s_1)\big)^{\omega- \gamma}$, we can use this estimate in the Lemma \ref{LemmaB2} and in the left hand side of inequality (\ref{EstimationB2}). We obtain in the same manner and with the same constants:
$$\frac{d}{ds_1}\|\psi(\cdot,s_1)\|_{L^\infty}\leq -C\big(r+K(s_0+s_1)\big)^{-\frac{(\omega- \gamma)}{n+\omega}}\|\psi(\cdot,s_1)\|_{L^\infty}^{1+\frac{1}{n+\omega}}.$$
We remark now to conclude that if  $\|\psi(\cdot,s_1)\|_{L^\infty}\leq (r+K(s_0+s_1))^{-(n+\gamma)}$ then the previous inequality is satisfied.\hfill $\blacksquare$\\

The crucial part of the proof of the Theorem \ref{Generalisacion} is given by the next proposition which gives us a control on the $L^1$-norm for a time $s_0+s_1$.
\begin{Proposition}[Second $L^1$-norm estimate]\label{Propo6}
Under the hypothesis of the Theorem (\ref{Generalisacion}) we have
$$\|\psi(\cdot,s_1)\|_{L^1}\leq \frac{v_n}{\big(r+K(s_0+s_1)\big)^{\gamma}}$$
\end{Proposition}
\textit{\textbf{Proof.}} Once we have at our disposal the concentration condition (\ref{Concentration2}) and the height condition (\ref{Linftyevolutionnotsmalltime}), the proof of the $L^1$ bound follows the same lines of the Corollary \ref{CorollaireL1}.\hfill $\blacksquare$
\subsection{The iteration}\label{SecIterationMol}
In sections \ref{SecEvolMol1} and \ref{SecEvolMol2} we studied respectively the evolution of small molecules from time $0$ to a small time $s_0$ and from this time $s_0$ to a larger time $s_0+s_1$ and we obtained a good $L^1$ control for such molecules. It is now possible to reapply the Theorem \ref{Generalisacion} in order to obtain a larger time control of the $L^1$ norm. The calculus of the $N$-th iteration will be exactly the same, although it will be necessary to make some modifications. 
\begin{Theoreme}\label{TheoIterationFin}
Assume that $\psi(x, s_N)$ is a solution of the system
\begin{equation*}
\left\lbrace
\begin{array}{rl}
\partial_{s_N} \psi(x,s_N)=& -\nabla\cdot(v\, \psi)(x,s_N)-\Lambda\psi(x,s_N)\\[5mm]
\psi(x,0)=& \psi(x,s_{N-1})  \qquad \qquad \mbox{with } s_{N-1}>0\\[5mm]
div(v)=&0 \quad \mbox{and }\; v\in L^{\infty}([0,T];bmo(\mathbb{R}^n))\quad \mbox{with } \underset{s_N\in [s_{N-1},T]}{\sup}\; \|v(\cdot,s_N)\|_{bmo}\leq \mu
\end{array}
\right.
\end{equation*}
If $\psi(x,s_{N-1})$ satisfies the three following conditions
\begin{eqnarray}
\int_{\mathbb{R}^n}|\psi(x,s_{N-1})||x-x(s_{N-1})|^\omega dx &\leq &(r+K(s_0+\cdots+s_{N-1}))^{\omega-\gamma}\nonumber\\
\|\psi(\cdot, s_{N-1})\|_{L^\infty}&\leq & \frac{1}{\left(r+ K(s_0+\cdots+s_{N-1}\right)^{n+\gamma}}\nonumber\\
\|\psi(\cdot, s_{N-1})\|_{L^1} &\leq & \frac{v_n}{\big(r+K(s_0+\cdots+ s_{N-1}\big)^{\gamma}}\nonumber
\end{eqnarray}
Then for all $0< s_N\leq\epsilon r$ small, we have the following estimates
\begin{eqnarray}
\int_{\mathbb{R}^n}|\psi(x,s_N)||x-x(s_N)|^\omega dx &\leq &(r+ K(s_0+\cdots+s_{N-1}+ s_N))^{\omega-\gamma}  \label{Concentration2Fin}\\
\|\psi(\cdot,s_N)\|_{L^\infty}&\leq & \frac{1}{\left(r+K(s_0+\cdots+s_{N-1}+s_N)\right)^{n+\gamma}}\label{LinftyevolutionnotsmalltimeFin}\\
\|\psi(\cdot,s_N)\|_{L^1} &\leq & \frac{v_n}{\big(r+K(s_0+\cdots+s_{N-1}+s_N)\big)^{\gamma}} \label{L1evolutionsmalltimeFin}
\end{eqnarray}
\end{Theoreme}
\begin{Remarque}
\emph{
\begin{itemize} 
\item[1)] Again, it is enough to assume that $(r+ K(s_0+\cdots+s_{N}))<1$, otherwise there is nothing to prove.
\item[2)] The new molecule's center $x(s_N)$ used in formula (\ref{Concentration2Fin}) is fixed by 
\begin{equation}\label{Defpointx_sfin}
\left\lbrace
\begin{array}{rl}
x'(s_N)=& \overline{v}_{B_{f_N}}=\frac{1}{|B_{f_N}|}\displaystyle{\int_{B_{f_N}}}v(y,s_N)dy\\[5mm]
x(0)=& x(s_{N-1}).
\end{array}
\right.
\end{equation}  
And here we noted $B_{f_N}=B(x(s_N),f_N)$ with $f_N$ a real valued function given by 
\begin{equation}\label{DefiFunctionFfin}
f_N= \big(r+K(s_0+\cdots+s_{N-1})\big).
\end{equation}
\end{itemize}
}
\end{Remarque}
\textit{\textbf{Proof of the Theorem \ref{TheoIterationFin}.}} We start with inequality (\ref{Concentration2Fin}). Let us write again $\Omega(x-x(s_N))=|x-x(s_N)|^{\omega}$ and $\psi(x)=\psi_+(x)-\psi_-(x)$, then we have $|\psi(x,s_N)|=|\psi_+(x,s_N)-\psi_-(x,s_N)|\leq \psi_+(x,s_N)+\psi_-(x,s_N)$ and we write
$$\int_{\mathbb{R}^n}|\psi(x,s_N)|\Omega(x-x(s_N))dx\leq\int_{\mathbb{R}^n}\psi_+(x,s_N)\Omega(x-x(s_N))dx+\int_{\mathbb{R}^n}\psi_-(x,s_N)\Omega(x-x(s_N))dx.$$
We only treat one of the integrals on the right-hand side above:
\begin{eqnarray*}
I&=&\left|\partial_{s_N} \int_{\mathbb{R}^n}\Omega(x-x(s_N))\psi_+(x,s_N)dx\right|.
\end{eqnarray*}
Following the same steps as before and using the fact that $v$ is divergence free, we obtain
\begin{equation*}
I=\left|\int_{\mathbb{R}^n}\nabla\Omega(x-x(s_N))\cdot(v-x'(s_N))\psi_+(x,s_N)-\Omega(x-x(s_N))\Lambda\psi_+(x,s_N)dx\right|.
\end{equation*}
Using the definition of $x'(s_N)$ given in (\ref{Defpointx_sfin}) and replacing $\Omega(x-x(s_N))$ by $|x-x(s_N)|^{\omega}$ we have
\begin{equation}\label{Estrellafin}
I\leq c \underbrace{\int_{\mathbb{R}^n}|x-x(s_N)|^{\omega-1}|v-\overline{v}_{B_{f_N}}| |\psi_+(x,s_N)|dx}_{I_1} + c\underbrace{\int_{\mathbb{R}^n}|x-x(s_N)|^{\omega-1}|\psi_+(x,s_N)|dx}_{I_2}.
\end{equation}
We will study separately each of the integrals $I_1$ and $I_2$ in the next lemmas:
\begin{Lemme}\label{Lemme41fin} For integral $I_1$ we have the estimate $I_1\leq C \mu\big(r+K(s_0+\cdots+s_{N-1})\big)^{\omega-\gamma-1}$.
\end{Lemme}
\textit{\textbf{Proof}.} We set $\mathbb{R}^n=B_{f_N}\cup\bigcup_{k\geq 1}E_k$ with
\begin{eqnarray}\label{Decoupagefin}
B_{f_N}&=& \{x\in \mathbb{R}^n: |x-x(s_N)|\leq f_N\},\\[5mm]
E_k&=& \{x\in \mathbb{R}^n: f_N2^{k-1}<|x-x(s_N)|\leq f_N 2^{k}\} \quad \mbox{for } k\geq 1.\nonumber
\end{eqnarray}
\begin{enumerate}
\item[(i)] \underline{Estimations over the ball $B_{f_N}$}. Applying Hölder inequality on integral $I_1$ we obtain
\begin{eqnarray*}
I_{1,B_{f_N}}=\int_{B_{f_N}}|x-x(s_N)|^{\omega-1}|v-\overline{v}_{B_{f_N}}| |\psi_+(x,s_N)|dx & \leq &\underbrace{\||x-x(s_N)|^{\omega-1}\|_{L^p(B_{f_N})}}_{(1)} \\
& \times &\underbrace{\|v-\overline{v}_{B_{f_N}}\|_{L^z(B_{f_N})}}_{(2)}\underbrace{\|\psi_+(\cdot, s_N)\|_{L^q(B_{f_N})}}_{(3)}
\end{eqnarray*}
where $\frac{1}{p}+\frac{1}{z}+\frac{1}{q}=1$ and $p,z,q> 1$. 
\begin{enumerate}
\item[$\bullet$] Observe that for $1<p<n/(1-\omega)$ we have $\||x-x(s_N)|^{\omega-1}\|_{L^p(B_{f_N})}\leq C f_N^{n/p+\omega-1}$.

\item[$\bullet$] By hypothesis we have $v(\cdot, s_N)\in bmo$, thus $\|v-\overline{v}_{B_{f_N}}\|_{L^z(B_{f_N})}\leq C|B_{f_N}|^{1/z}\|v(\cdot,s_N)\|_{bmo}$, now since\\
 $\underset{s_N\in [s_{N-1},T]}{\sup}\; \|v(\cdot,s_N)\|_{bmo}\leq \mu$ we write $\|v-\overline{v}_{B_{f_N}}\|_{L^z(B_{f_N})}\leq C f_N^{n/z}\mu$.

\item[$\bullet$] Using the maximum principle for $L^q$ norms we have $\|\psi_+(\cdot, s_N)\|_{L^q(B_{f_N})}\leq  \|\psi(\cdot, s_{N-1})\|_{L^1}^{1/q}\|\psi(\cdot, s_{N-1})\|_{L^\infty}^{1-1/q}.$
\end{enumerate}
We gather all these inequalities in order to obtain the following estimation for $I_{1,B_{f_N}}$:
\begin{equation*}
I_{1,B_{f_N}}\leq C\mu f_N^{n(1-1/q)+\omega-1}\|\psi(\cdot, s_{N-1})\|_{L^1}^{1/q}\|\psi(\cdot, s_{N-1})\|_{L^\infty}^{1-1/q}.
\end{equation*}
\item[(ii)] \underline{Estimations for the dyadic corona $E_k$}. Let us note $I_{1,E_k}$ the integral
$$I_{1,E_k}=\int_{E_k}|x-x(s_N)|^{\omega-1}|v-\overline{v}_{B_{f_N}}| |\psi_+(x,s_N)|dx.$$
Since over $E_k$ we have $|x-x(s_N)|^{\omega-1}\leq C 2^{k(\omega-1)}f_N^{\omega-1}$ we write
\begin{eqnarray*}
I_{1,E_k}&\leq & C2^{k(\omega-1)}f_N^{\omega-1}\left(\int_{E_k}|v-\overline{v}_{B(f_N2^k)}| |\psi_+(x,s_N)|dx+\int_{E_k}|\overline{v}_{B_{f_N}}-\overline{v}_{B(f_N2^k)}| |\psi_+(x,s_N)|dx\right)\\
&\leq& C2^{k(\omega-1)}f_N^{\omega-1}\left(\int_{B(f_N2^k)}|v-\overline{v}_{B(f_N2^k)}| |\psi_+(x,s_N)|dx+\int_{B(f_N2^k)}|\overline{v}_{B_{f_N}}-\overline{v}_{B(f_N2^k)}| |\psi_+(x,s_N)|dx\right),
\end{eqnarray*}
with  $B(f_N2^k)=B(x(s_N), f_N2^k)$. Now, since $v(\cdot, s_N)\in bmo$, using Proposition \ref{PropoBMO1} we have $|\overline{v}_{B_{f_N}}-\overline{v}_{B(f_N2^k)}| \leq Ck\|v(\cdot,s_N)\|_{bmo}\leq Ck\mu$ and we can write
\begin{eqnarray*}
I_{1,E_k}&\leq &C2^{k(\omega-1)}f_N^{\omega-1}\left(\int_{B(f_N2^k)}|v-\overline{v}_{B(f_N2^k)}| |\psi_+(x,s_N)|dx+Ck\mu\|\psi_+(\cdot,s_N)\|_{L^1}\right)\\[5mm]
&\leq & C2^{k(\omega-1)}f_N^{\omega-1}\left(\|\psi_+(\cdot,s_N)\|_{L^{a_0}} \|v-\overline{v}_{B(f_N2^k)}\|_{L^{\frac{a_0}{a_0-1}}} +Ck\mu\;\|\psi_+(\cdot,s_{N-1})\|_{L^1}\right)
\end{eqnarray*}
where we used Hölder inequality with $1<a_0<\frac{n}{n+(\omega-1)}$ and maximum principle for the last term above. Using again the properties of $bmo$ spaces we have
$$I_{1,E_k}\leq  C2^{k(\omega-1)}f_N^{\omega-1}\left(\|\psi_+(\cdot,s_{N-1})\|_{L^1}^{1/a_0}\|\psi_+(\cdot,s_{N-1})\|_{L^\infty}^{1-1/a_0} |B(f_N2^k)|^{1-1/a_0}\|v(\cdot,s_N)\|_{bmo} +Ck\mu\|\psi(\cdot,s_{N-1})\|_{L^1}\right).$$
Since $\|v(\cdot,s_{N})\|_{bmo}\leq \mu$ and since $1<a_0<\frac{n}{n+(\omega-1)}$, we have $n(1-1/a_0)+(\omega-1)<0$, so that, summing over each dyadic corona $E_k$, we obtain
\begin{equation*}
\sum_{k\geq 1}I_{1,E_k}\leq C \mu\left(f_N^{n(1-1/a_0)+\omega-1}\|\psi(\cdot,s_{N-1})\|_{L^1}^{1/a_0}\|\psi(\cdot,s_{N-1})\|_{L^\infty}^{1-1/a_0}+f_N^{\omega-1}\|\psi(\cdot,s_{N-1})\|_{L^1}\right).
\end{equation*}
\end{enumerate}
We finally obtain the following inequalities:
\begin{eqnarray}
I_1&=& I_{1,B_{f_N}}+\sum_{k\geq 1}I_{1,E_k}\label{FinalEstimateI_1fin}\\
&\leq &C\mu \underbrace{f_N^{n(1-1/q)+\omega-1}\|\psi(\cdot, s_{N-1})\|_{L^1}^{1/q}\|\psi(\cdot, s_{N-1})\|_{L^\infty}^{1-1/q}}_{(a)}\nonumber\\
& &+ C \mu \left(\underbrace{f_N^{n(1-1/a_0)+\omega-1}\|\psi(\cdot,s_{N-1})\|_{L^1}^{1/a_0}\|\psi(\cdot,s_{N-1})\|_{L^\infty}^{1-1/a_0}}_{(b)} +\underbrace{f_N^{\omega-1}\|\psi(\cdot,s_{N-1})\|_{L^1}}_{(c)}\right)\nonumber
\end{eqnarray}
Each of the terms $(a)$, $(b)$ and $(c)$ above is bounded by the quantity $\big(r+K(s_0+\cdots+s_{N-1})\big)^{\omega-\gamma-1}$:
\begin{itemize}
\item for the first term (a) by the hypothesis on the initial data $\psi(\cdot,s_{N-1})$ and the definition of $f_N$ given in (\ref{DefiFunctionFfin}) we have:
\begin{eqnarray*}
f_N^{n(1-1/q)+\omega-1}\|\psi(\cdot, s_0)\|_{L^1}^{1/q}\|\psi(\cdot, s_0)\|_{L^\infty}^{1-1/q}& \leq & \big(r+K(s_0+\cdots+s_{N-1})\big)^{[n(1-1/q)+\omega-1]-\frac{\gamma}{q}-(n+ \gamma)(1-1/q)}\\
& \leq & \big(r+K(s_0+\cdots+s_{N-1})\big)^{\omega-\gamma-1}.
\end{eqnarray*}
\item For the second term (b) we have:
\begin{eqnarray*}
f_N^{n(1-1/a_0)+\omega-1}\|\psi(\cdot,s_0)\|_{L^1}^{1/a_0}\|\psi(\cdot,s_0)\|_{L^\infty}^{1-1/a_0} & \leq & \big(r+K(s_0+\cdots+s_{N-1})\big)^{[n(1-1/a_0)+\omega-1]-\frac{\gamma}{a_0}-(n+ \gamma)(1-1/a_0)}\\
& \leq & \big(r+K(s_0+\cdots+s_{N-1})\big)^{\omega-\gamma-1}.
\end{eqnarray*}
\item Finally, for the last term (c) we write
\begin{eqnarray*}
f_N^{\omega-1}\|\psi(\cdot,s_0)\|_{L^1}&\leq &f_N^{\omega-1}(r+K(s_0+\cdots+s_{N-1}))^{-\gamma}= (r+K(s_0+\cdots+s_{N-1}))^{\omega-\gamma-1}\\
\end{eqnarray*}
\end{itemize}
Gathering these estimates on $(a), (b)$ and $(c)$, and getting back to (\ref{FinalEstimateI_1fin}) we finally obtain 
$$I_1\leq C\mu \big(r+K(s_0+\cdots+s_{N-1})\big)^{\omega-\gamma-1}.$$ 
The Lemma \ref{Lemme41fin} is proven. \hfill$\blacksquare$\\
\begin{Lemme}\label{Lemme42fin}
For integral $I_2$ in inequality (\ref{Estrellafin}) we have the following estimate
\begin{equation*}
I_2=\int_{\mathbb{R}^n}|x-x(s_N)|^{\omega-1}|\psi_+(x,s_N)|dx\leq C\big(r+K(s_0+\cdots+s_{N-1})\big)^{\omega-\gamma-1}.
\end{equation*}
\end{Lemme}
\textbf{\textit{Proof.}} As for the Lemma \ref{Lemme41fin}, we consider $\mathbb{R}^n$ as the union of a ball with dyadic coronas centered on $x(s_N)$ (cf. (\ref{Decoupagefin})). \\
\begin{enumerate}
\item[(i)] \underline{Estimations over the ball $B_{f_N}$}. Applying Hölder inequality with $1< a_1<n/(1-\omega)$ and maximum principle we have
\begin{eqnarray*}
I_{2,B_{f_N}}=\int_{B_{f_N}}|x-x(s_N)|^{\omega-1}|\psi_+(x,s_N)|dx &\leq &\||x-x(s_N)|^{\omega-1}\|_{L^{a_1}(B_{f_N})}\|\psi_+(\cdot, s_N)\|_{L^{b_1}(B_{f_N})}\\
&\leq & Cf_N^{n/a_1+\omega-1}\|\psi_+(\cdot, s_{N-1})\|_{L^1}^{1-1/a_1}\|\psi_+(\cdot, s_{N-1})\|_{L^\infty}^{1/a_1}.
\end{eqnarray*}
Thus, we can write:
\begin{equation}\label{Bola2fin}
I_{2,B_{f_N}}\leq Cf_N^{n/a_1+\omega-1}\|\psi_+(\cdot, s_{N-1})\|_{L^1}^{1-1/a_1}\|\psi_+(\cdot, s_{N-1})\|_{L^\infty}^{1/a_1}.
\end{equation}
\item[(ii)]  \underline{Estimations for the dyadic corona $E_k$}. Here we have
\begin{eqnarray*}
I_{2,E_k}=\int_{E_k}|x-x(s_N)|^{\omega-1}|\psi_+(x,s_N)|dx & \leq & C2^{k(\omega-1)}f_N^{\omega-1}\int_{E_k}|\psi_+(x,s_N)|dx\leq C2^{k(\omega-1)}f_N^{\omega-1}\|\psi_+(\cdot, s_N)\|_{L^1}\\
&\leq & C2^{k(\omega-1)}f_N^{\omega-1}\|\psi(\cdot, s_{N-1})\|_{L^1}
\end{eqnarray*}
Since $0<\gamma<\omega<1$ we have $\omega-1<0$ and summing over $k\geq 1$ we obtain
\begin{equation*}
\sum_{k\geq 1}I_{2,E_k}=\sum_{k\geq 1}\int_{E_k}|x-x(s_N)|^{\omega-1}|\psi_+(x,s_N)|dx \leq Cf_N^{\omega-1}\|\psi(\cdot, s_{N-1})\|_{L^1}.
\end{equation*}
We finally obtain
\begin{equation}\label{Coronak2fin}
\sum_{k\geq 1}I_{2,E_k}\leq C f_N^{\omega-1}\|\psi(\cdot, s_{N-1})\|_{L^1}.
\end{equation}
\end{enumerate}
To finish the proof of the Lemma \ref{Lemme42fin}  we glue together (\ref{Bola2fin}) and (\ref{Coronak2fin}) and we obtain
$$I_2=I_{2,B_{f_N}}+\sum_{k\geq 1}I_{2,E_k}\leq C\left(\underbrace{f_N^{n/a_1+\omega-1}\|\psi_+(\cdot, s_{N-1})\|_{L^1}^{1-1/a_1}\|\psi_+(\cdot, s_{N-1})\|_{L^\infty}^{1/a_1}}_{(d)}+\underbrace{f_N^{\omega-1}\|\psi(\cdot, s_{N-1})\|_{L^1}}_{(e)}\right)$$
\begin{itemize}
\item For the term $(d)$ we write
\begin{eqnarray*}
f_N^{n/a_1+\omega-1}\|\psi(\cdot, s_{N-1})\|_{L^1}^{1-1/a_1}\|\psi(\cdot, s_{N-1})\|_{L^\infty}^{1/a_1}& \leq & \big(r+K(s_0+\cdots+s_{N-1})\big)^{[n/a_1+\omega-1]-\gamma(1-1/a_1)-\frac{n+\gamma}{a_1}}\\
& \leq & \big(r+K(s_0+\cdots+s_{N-1})\big)^{\omega-\gamma-1}
\end{eqnarray*}
\item To treat the term $(e)$ it is enough to apply the same arguments used to prove the part $(c)$ above. 
\end{itemize}
Finally, we obtain $I_2=I_{2,B_{f_N}}+\displaystyle{\sum_{k\geq 1}}I_{2, E_k}\leq  C\big(r+K(s_0+\cdots+s_{N-1})\big)^{\omega-\gamma-1}$ and the Lemma \ref{Lemme42fin} is proven.\hfill$\blacksquare$\\

Now, using Lemmas \ref{Lemme41fin} and \ref{Lemme42fin} and getting back to the estimate (\ref{Estrellafin}) we have
\begin{equation}\label{FinalEstimatefin}
\left|\partial_{s_N} \int_{\mathbb{R}^n}\Omega(x-x(s_N))\psi_+(x,s_N)dx\right| \leq  C \left(\mu+1\right)\big(r+K(s_0+\cdots+s_{N-1})\big)^{\omega-\gamma-1}
\end{equation}

This estimation is compatible with the estimate (\ref{Concentration2Fin}) for $0\leq s_N\leq \epsilon r$ small:  write
$\phi=\left(r+K(s_0+\cdots+s_{N})\right)^{\omega-\gamma}$ and linearize this expression with respect to $s_N$:
$$\phi \thickapprox (r+K(s_0+\cdots+s_{N-1}))^{\omega-\gamma}\left(1+K(\omega-\gamma)\frac{s_N}{(r+K(s_0+\cdots+s_{N-1}))}\right).$$
Taking the derivative of $\phi$ with respect to $s_N$ we have $\phi' \thickapprox K(\omega-\gamma) \big(r+K(s_0+\cdots+s_{N-1})\big)^{\omega-\gamma-1}$ and with the condition (\ref{SmallConstants}) on $K(\omega-\gamma)$ we obtain that (\ref{FinalEstimatefin}) is bounded by $\phi'$ and the concentration condition follows.\\

Let us study now the height condition. As long as this concentration condition is bounded by $(r+K(s_0+\cdots+s_{N}))^{\omega-\gamma}$, we can deduce from it the $L^\infty$ estimate
$$\|\psi(\cdot, s_N)\|_{L^\infty}\leq \frac{1}{\big(r+K(s_0+\cdots+s_{N})\big)^{\gamma}}.$$
Indeed, following the same steps of the Proposition \ref{Propo2} and Lemma \ref{LemmaB2}, we obtain the inequality:
$$\frac{d}{ds_N}\|\psi(\cdot,s_N)\|_{L^\infty}\leq -C\big(r+K(s_0+\cdots+s_N)\big)^{-\frac{(\omega- \gamma)}{n+\omega}}\|\psi(\cdot,s_N)\|_{L^\infty}^{1+\frac{1}{n+\omega}}.$$
It suffices to note that if $\|\psi(\cdot,s_N)\|_{L^\infty}\leq \frac{1}{\big(r+K(s_0+\cdots+s_N)\big)^{n+\gamma}}$, then $\|\psi(\cdot,s_N)\|_{L^\infty}$ satifies this inequality.\\

Finally, proceeding as in the corollary \ref{CorollaireL1} with these two inequalities we have a $L^1$ estimate for time $s_0+\cdots+s_{N}$:
$$\|\psi(\cdot,s_N)\|_{L^1}\leq \frac{v_n}{\big(r+K(s_0+\cdots+s_{N})\big)^{\gamma}}.$$
\hfill$\blacksquare$\\

\textit{\textbf{End of the proof of the Theorem \ref{TheoL1estimate}.}} We have proved that is possible to control the $L^1$ behavior of the molecules from $0$ to a time $s_0$, from time $s_0$ to time $s_1$, and by iteration from time $s_{N-1}$ to time $s_N$. Observe now that the smallness of $r$ and of the times $s_0,...,s_N$ can be compensated by the number of iterations $N$ in the following sense: fix a small $0<r<1$ and iterate as explained before. Since each small time $s_0,...,s_N $ is of order $\epsilon r$, we have $s_0+\cdots+s_{N}\sim N \epsilon r$. Thus, we will stop the iterations as soon as $ N r\geq T_0$. Of course, the number of iterations $N=N(r)$ will depend on the smallness of the molecule's size $r$, more specifically it is enough to set $N(r)\sim \frac{T_0}{r}$. Proceeding this way we will obtain $\|\psi(\cdot,s_N)\|_{L^1}\leq C T_0^{-\gamma}<+\infty$. Note in particular that, once this estimate is available, for bigger times it is enough to apply the maximum principle.\\

Finally, and for all $r>0$, we obtain after a time $T_0$ a $L^1$ control for small molecules and we finish the proof of the Theorem \ref{TheoL1estimate}. \hfill$\blacksquare$
\section{Existence and uniqueness for $L^p$ initial data}\label{SecExiUnic}
In this section we will study existence and uniqueness for weak solution of equation (\ref{QG}) with initial data $\theta_0\in L^p(\mathbb{R}^n)$ where $p\geq 2$. Remark that equation (\ref{QG}) differs mainly from the backward equation (\ref{Evolution1}) by the sign of the velocity. Since the velocity $v$ is a data for the problem, it is equivalent to consider $-v$ instead of $v$, thus, for simplicity, we fix velocity's sign in the next way:
\begin{equation}\label{QG+}
\left\lbrace
\begin{array}{l}
\partial_t \theta(x,t)+ \nabla\cdot (v\;\theta)(x,t)+\Lambda\theta(x,t)=0\\[5mm]
\theta(x,0)=\theta_0(x)\in L^p(\mathbb{R}^n)\\[5mm]
div(v)=0 \quad \mbox{ and } v\in L^{\infty}([0,T];bmo(\mathbb{R}^n)).
\end{array}
\right.
\end{equation}
\subsection{Viscosity Solutions}\label{SecViscos}
The term \textit{Viscosity Solutions} is taken from \cite{Cordoba} and it refers to weak solutions of (\ref{QG+}) which are the weak limit, as $\varepsilon\longrightarrow 0$, of a sequence of solutions of problems
\begin{equation}\label{SistApprox}
\left\lbrace
\begin{array}{l}
\partial_t \theta(x,t)+ \nabla\cdot(v_{\varepsilon}\;\theta)(x,t)+\Lambda\theta(x,t)=\varepsilon \Delta \theta(x,t)\\[5mm]
\theta(x,0)=\theta_0(x)\\[5mm]
div(v)=0 \quad \mbox{ and } v\in L^{\infty}([0,T];  L^{\infty}(\mathbb{R}^n)).
\end{array}
\right.
\end{equation}
where $v_{\varepsilon}$ is defined by  $v_{\varepsilon}=v\ast \omega_{\varepsilon}$ with $\omega_{\varepsilon}(x)=\varepsilon^{-n}\omega(x/\varepsilon)$ and $\omega\in \mathcal{C}^{\infty}_0(\mathbb{R}^n)$ is a function such that $\displaystyle{\int_{\mathbb{R}^n}}\omega(x)dx=1$.
\begin{Remarque}\label{Rem_Approx}
\emph{Observe that we fixed here the velocity $v$ such that $v\in L^{\infty}([0,T];  L^{\infty}(\mathbb{R}^n))$. This is not very restrictive because by Proposition \ref{TheoApproxbmo} we can construct a sequence $v_k\in L^\infty$ that converge weakly to $v$ in $bmo$. }
\end{Remarque}
Problem (\ref{SistApprox}) admits the following equivalent integral representation:
\begin{equation}\label{FormIntegr}
\theta(x,t)=e^{\varepsilon t\Delta}\theta_0(x)-\int_{0}^t e^{\varepsilon (t-s)\Delta}\nabla \cdot(v_\varepsilon\; \theta)(x,s)ds-\int_{0}^t e^{\varepsilon (t-s)\Delta}\Lambda \theta(x,s)ds
\end{equation}
For a proof of this fact see \cite{Marchand} or \cite{PGLR}. We will use then the Banach contraction scheme and for this we will consider the space $L^{\infty}([0,T]; L^{p}(\mathbb{R}^n))$ with the norm $\|f\|_{L^\infty (L^p)}=\underset{t\in [0,T]}{\sup}\|f(\cdot, t)\|_{L^p}$.
\begin{Theoreme}[Local existence]\label{TheoPointFixe}
Let $2\leq p<+\infty$ and let $\theta_0$ and $v$ be two functions such that $\theta_0\in L^p(\mathbb{R}^n)$, $div(v)=0$ and $v\in L^{\infty}([0,T']; L^{\infty}(\mathbb{R}^n))$. If initial data satisfies $\|\theta_0\|_{L^p}\leq K$ and if $T'$ is a time small enough such that
$$C\left(\frac{T'^{1/2}}{\varepsilon^{1/2}}\|v\|_{L^\infty (L^\infty)} + \frac{T'^{1/2}}{\varepsilon^{1/2}}\right)\leq 1/2,$$ 
then (\ref{FormIntegr}) has a unique solution $\theta \in L^{\infty}([0,T']; L^{p}(\mathbb{R}^n))$ on the closed ball $\overline{B}(0,2K)\subset L^{\infty}([0,T']; L^{p}(\mathbb{R}^n))$. 
\end{Theoreme}
\textit{\textbf{Proof}.}
We note $L_\varepsilon(\theta)$ and  $N^{v}_\varepsilon(\theta)$ the quantities
\begin{eqnarray*}
L_\varepsilon(\theta)(x,t)= \int_{0}^t e^{\varepsilon (t-s)\Delta}\Lambda \theta(x,s)ds & \mbox{ and } & 
N^{v}_\varepsilon(\theta)(x,t)=\int_{0}^t e^{\varepsilon (t-s)\Delta}\nabla \cdot(v_\varepsilon\; \theta)(x,s)ds. 
\end{eqnarray*}
\begin{Lemme} If $f\in L^{\infty}([0,T']; L^{p}(\mathbb{R}^n))$, then
\begin{equation}\label{Maj2}
\|L_\varepsilon(f)\|_{L^\infty (L^p)}\leq C \frac{T'^{1/2}}{\varepsilon^{1/2}}\; \|f\|_{L^\infty (L^p)}
\end{equation}
\end{Lemme}
\textit{\textbf{Proof.}}
We write
$$\|L_\varepsilon(f)\|_{L^\infty (L^p)}=\underset{0<t<T'}{\sup} \left\|\int_{0}^t e^{\varepsilon (t-s)\Delta}\Lambda f(\cdot,s)ds\right\|_{L^p}=\underset{0<t<T'}{\sup} \left\|\int_{0}^t \Lambda f\ast h_{\varepsilon (t-s)}(\cdot,s)ds\right\|_{L^p}$$
where we noted $h_t$ the heat kernel on $\mathbb{R}^n$. Then we have the estimates
\begin{eqnarray*}
\|L_\varepsilon(f)\|_{L^\infty (L^p)}&\leq &\underset{0<t<T'}{\sup} \int_{0}^t \left\|f(\cdot, s)\right\|_{L^p}  \left\|\Lambda h_{\varepsilon (t-s)}\right\|_{L^1}ds\leq \left\|f\right\|_{L^\infty (L^p)} \underset{0<t<T'}{\sup} \int_{0}^t  C(\varepsilon(t-s))^{-1/2}ds\leq  C \frac{T'^{1/2}}{\varepsilon^{1/2}}\; \|f\|_{L^\infty (L^p)}.\quad \blacksquare
\end{eqnarray*}
\begin{Lemme} If $f\in L^{\infty}([0,T']; L^{p}(\mathbb{R}^n))$ and if $v\in  L^{\infty}([0,T']; L^{\infty}(\mathbb{R}^n)) $, then 
\begin{equation}\label{Maj3}
\|N^v_\varepsilon(f)\|_{L^\infty (L^p)}\leq C \sqrt{\frac{T'}{\varepsilon}}\;\|v\|_{L^\infty (L^\infty)} \|f\|_{L^\infty (L^p)}
\end{equation}
\end{Lemme}
\textit{\textbf{Proof.}} We write:
\begin{eqnarray*}
\|N^v_\varepsilon(f)\|_{L^\infty (L^p)}& =&\underset{0<t<T'}{\sup} \left\|\int_{0}^t e^{\varepsilon (t-s)\Delta} \nabla \cdot(v_\varepsilon f)(\cdot,s)ds\right\|_{L^p}=\underset{0<t<T'}{\sup} \left\|\int_{0}^t \nabla \cdot(v_\varepsilon f)\ast h_{\varepsilon (t-s)}(\cdot,s)ds\right\|_{L^p}\\
&\leq & \underset{0<t<T'}{\sup}\int_{0}^t  \left\|v_\varepsilon f(\cdot,s)\right\|_{L^p} \left\|\nabla h_{\varepsilon (t-s)}\right\|_{L^1} ds\leq  \underset{0<t<T'}{\sup}\int_{0}^t  \left\|v_\varepsilon(\cdot,s)\right\|_{L^\infty}  \left\|f(\cdot,s)\right\|_{L^p} C(\varepsilon(t-s))^{-1/2} ds\\
&\leq &  \left\|f\right\|_{L^\infty (L^p)} \left\|v\right\|_{L^\infty (L^\infty)} \underset{0<t<T'}{\sup}\int_{0}^t C(\varepsilon(t-s))^{-1/2} ds\leq  \left\|f\right\|_{L^\infty (L^p)} \|v\|_{L^\infty (L^\infty)} C \sqrt{\frac{T'}{\varepsilon}}.\qquad \blacksquare
\end{eqnarray*}
Finally, since $e^{\varepsilon t \Delta}$ is a contraction operator, estimate $\|e^{\varepsilon t \Delta}f\|_{L^p}\leq \|f\|_{L^p}$ is valid for all function $f\in L^{p}(\mathbb{R}^n)$ with $1\leq p\leq +\infty$, for all $t>0$ and all $\varepsilon>0$. Thus, we have
\begin{equation}\label{Maj1}
\|e^{\varepsilon t \Delta}f\|_{L^\infty (L^p)}\leq \|f\|_{L^p}.
\end{equation}
To apply the Banach contraction scheme, let us now construct a sequence of functions in the following way
$$\theta_{n+1}(x,t)=e^{\varepsilon t\Delta}\theta_0(x)-L_{\varepsilon}(\theta_n)(x,t)-N^v_{\varepsilon}(\theta_n)(x,t)$$
and we take the $L^\infty (L^p)$-norm of this expression to obtain
$$ \|\theta_{n+1}\|_{L^\infty (L^p)}\leq\|e^{\varepsilon t\Delta}\theta_0\|_{L^\infty (L^p)}+\|L_{\varepsilon}(\theta_n)\|_{L^\infty (L^p)}+\|N^v_{\varepsilon}(\theta_n)\|_{L^\infty (L^p)}$$
Using estimates (\ref{Maj2}), (\ref{Maj3}) and (\ref{Maj1}) we have
$$ \|\theta_{n+1}\|_{L^\infty (L^p)}\leq \|\theta_0\|_{L^p}+C\left(\frac{T'^{1/2}}{\varepsilon^{1/2}}\|v\|_{L^\infty (L^\infty)} + \frac{T'^{1/2}}{\varepsilon^{1/2}}\right)\|\theta_n\|_{L^\infty (L^p)}$$
Thus, if $\|\theta_0\|_{L^p}\leq K$ and with the definition of $T'$, we have by iteration that $\|\theta_{n+1}\|_{L^\infty (L^p)}\leq 2 K$: the sequence $(\theta_n)_{n\in \mathbb{N}}$ constructed from initial data $\theta_0$ belongs to the closed ball $\overline{B}(0, 2K)$. In order to finish this proof, let us show that $\theta_n \longrightarrow \theta$ in $L^{\infty}([0,T']; L^{p}(\mathbb{R}^n))$. For this we write
$$\|\theta_{n+1}-\theta_n\|_{L^\infty (L^p)}\leq \|L_\varepsilon(\theta_{n}-\theta_{n-1})\|_{L^\infty (L^p)}+\|N^v_\varepsilon(\theta_{n}-\theta_{n-1})\|_{L^\infty (L^p)}$$
and using previous lemmas we have
$$\|\theta_{n+1}-\theta_n\|_{L^\infty (L^p)}\leq C\left(\frac{T'^{1/2}}{\varepsilon^{1/2}}\|v\|_{L^\infty (L^\infty)} + \frac{T'^{1/2}}{\varepsilon^{1/2}}\right)\|\theta_{n}-\theta_{n-1}\|_{L^\infty (L^p)}$$
so, by iteration we obtain
$$\|\theta_{n+1}-\theta_n\|_{L^\infty (L^p)}\leq \left[C\left(\frac{T'^{1/2}}{\varepsilon^{1/2}}\|v\|_{L^\infty (L^\infty)} + \frac{T'^{1/2}}{\varepsilon^{1/2}}\right)\right]^{n}\|\theta_{1}-\theta_0\|_{L^\infty (L^p)}$$
hence, with the definition of $T'$ it comes
$\|\theta_{n+1}-\theta_n\|_{L^\infty (L^p)}\leq \left(\frac{1}{2}\right)^{n}\|\theta_{1}-\theta_0\|_{L^\infty (L^p)}$.
Finally, if $n\longrightarrow +\infty$, the sequence $(\theta_n)_{n\in \mathbb{N}}$ convergences towards $\theta$ in $L^\infty([0,T'];L^p(\mathbb{R}^n))$. Since it is a Banach space we deduce uniqueness for the solution $\theta$ of problem (\ref{FormIntegr}).\hfill$\blacksquare$
\begin{Corollaire}\label{CorDepContinue}
The solution constructed above depends continuously on initial data $\theta_0$.
\end{Corollaire}
\textit{\textbf{Proof.}}
Let $\varphi_0\in L^{p}(\mathbb{R}^n)$ be an initial data and let $\varphi$ be the associated solution. We write
$$\theta(x,t)-\varphi(x,t)=e^{\varepsilon t\Delta}(\theta_0(x)-\varphi_0(x))-L_\varepsilon(\theta-\varphi)(x,t)-N^v_\varepsilon(\theta-\varphi)(x,t)$$
Taking $L^\infty (L^p)$-norm in formula above and applying the same previous calculations one obtains
\begin{equation}\label{ForUniBien}
\|\theta-\varphi\|_{L^\infty (L^p)}\leq  \|\theta_0-\varphi_0\|_{L^p}+ C_0\|\theta-\varphi\|_{L^\infty (L^p)}
\end{equation}
This shows continuous dependence of the solution since $C_0=\left(C\left(\frac{T'^{1/2}}{\varepsilon^{1/2}}\|v\|_{L^\infty (L^\infty)} + \frac{T'^{1/2}}{\varepsilon^{1/2}}\right)\right)\leq 1/2$.\hfill$\blacksquare$\\

Once we obtain a local result, global existence easily follows by a simple iteration since problems studied here (equations (\ref{QG}), (\ref{QG+}) or (\ref{SistApprox})) are linear as velocity $v$ does not depend on $\theta$. 
\begin{Remarque}
\emph{Solutions $\theta(\cdot,\cdot)$ constructed above depends on $\varepsilon$ and it will be more convenient to note them as $\theta^{(\varepsilon)}(\cdot,\cdot)$. For the time being, we will just note them $\theta(\cdot,\cdot)$.}
\end{Remarque}
We study now the regularity of solutions constructed by this method.
\begin{Theoreme} Solutions of the approximated problem (\ref{SistApprox}) are smooth.
\end{Theoreme}
\textit{\textbf{Proof}.}
By iteration we will prove that
$$\theta \in \bigcap_{0<T_0<T_1<t<T_2<T^\ast}L^\infty([0,t]; W^{\frac{k}{2},p}(\mathbb{R}^n))\quad \mbox{ for all } k\geq 0.$$ 
Remark that this is true for $k=0$. So let us assume that it is also true for $k>0$ and we will show that it is still true for $k+1$. Set $t$ such that $0<T_0<T_1<t<T_2<T^\ast$ and let us consider the next problem
$$\theta(x,t)=e^{\varepsilon (t-T_0)\Delta}\theta(x, T_0)-\int_{T_0}^t e^{\varepsilon (t-s)\Delta}\nabla \cdot(v_\varepsilon\; \theta)(x,s)ds-\int_{T_0}^t e^{\varepsilon (t-s)\Delta}\Lambda  \theta(x,s)ds
$$
We have then the following estimate
\begin{eqnarray*}
\|\theta\|_{L^\infty (W^{\frac{k+1}{2},p})}&\leq& \|e^{\varepsilon (t-T_0)\Delta}\theta(\cdot, T_0)\|_{L^\infty (W^{\frac{k+1}{2},p})}\\[5mm]
& &+\left\|\int_{T_0}^t e^{\varepsilon (t-s)\Delta}\nabla \cdot(v_\varepsilon\; \theta)(\cdot,s)ds\right\|_{L^\infty (W^{\frac{k+1}{2},p})}
+\left\|\int_{T_0}^t e^{\varepsilon (t-s)\Delta}\Lambda \theta(\cdot,s)ds\right\|_{L^\infty (W^{\frac{k+1}{2},p})}
\end{eqnarray*}
Now, we will treat separately each of the previous terms. 
\begin{enumerate}
\item[(i)] For the first one we have
$$\|e^{\varepsilon (t-T_0)\Delta}\theta(\cdot, T_0)\|_{W^{\frac{k+1}{2},p}}=\|\theta(\cdot, T_0)\ast \Lambda^{k+1}h_{\varepsilon (t-T_0)}\|_{L^p}\leq\|\theta(\cdot, T_0)\|_{L^p}\| \Lambda^{k+1}h_{\varepsilon (t-T_0)}\|_{L^1} $$ 
where $h_t$ is the heat kernel, so we can write
\begin{equation*}
\|e^{\varepsilon (t-T_0)\Delta}\theta(\cdot, T_0)\|_{L^\infty (W^{\frac{k+1}{2},p})}\leq C\|\theta(\cdot, T_0)\|_{L^p}\sup\left\{ \left[\varepsilon (t-T_0)\right]^{- \frac{k+1}{4}}; 1\right\}
\end{equation*}
\item[(ii)] For the second term, one has
\begin{eqnarray*}
I=\left\|\int_{T_0}^t e^{\varepsilon (t-s)\Delta}\nabla \cdot(v_\varepsilon\; \theta)(\cdot,s)ds\right\|_{W^{\frac{k+1}{2},p}}\leq \int_{T_0}^t\|\nabla  \cdot (v_{\varepsilon}\;\theta)\ast h_{\varepsilon(t-s)}\|_{W^{\frac{k+1}{2},p}}ds
\leq C\int_{T_0}^t\|v_\varepsilon\;\theta(\cdot,s)\|_{W^{\frac{k}{2},p}} \left[\varepsilon (t-s)\right]^{- \frac{3}{4}} ds
\end{eqnarray*}
Remark that we have here the estimations below for $N\geq k/2$
\begin{eqnarray*}
\|v_\varepsilon\theta(\cdot,s)\|_{W^{\frac{k}{2},p}}&\leq& \|v_\varepsilon(\cdot,s)\|_{\mathcal{C}^N} \|\theta(\cdot,s)\|_{W^{\frac{k}{2},p}}\leq C \varepsilon^{-N}\|v(\cdot,s)\|_{L^\infty}\|\theta(\cdot,s)\|_{W^{\frac{k}{2},p}}
\end{eqnarray*}
hence, we can write
\begin{eqnarray*}
I&\leq &C \|v\|_{L^\infty (L^\infty)} \|\theta\|_{L^\infty (W^{\frac{k}{2},p})}\int_{T_0}^t  \varepsilon^{-N}\sup\left\{  \left[\varepsilon (t-s)\right]^{- \frac{3}{4}} ;1\right\}ds
\end{eqnarray*}
\item[(iii)] Finally, for the last term we have
\begin{eqnarray*}
\left\|\int_{T_0}^t e^{\varepsilon (t-s)\Delta}\Lambda  \theta(\cdot,s)ds\right\|_{W^{\frac{k+1}{2},p}} &\leq & C\int_{T_0}^t \| \theta(\cdot,s) \|_{W^{\frac{k}{2},p}}\left[\varepsilon (t-s)\right]^{- \frac{3}{4}} ds
\end{eqnarray*}
\begin{eqnarray*}
\left\|\int_{T_0}^t e^{\varepsilon (t-s)\Delta}\Lambda  \theta(\cdot,s)ds\right\|_{L^\infty (W^{\frac{k+1}{2},p})}&\leq & C\|\theta \|_{L^{\infty}(W^{\frac{k}{2},p})}\int_{T_0}^t \sup\left\{  \left[\varepsilon (t-s)\right]^{- \frac{3}{4}} ;1\right\}ds.
\end{eqnarray*}
\end{enumerate}
Now, with formulas (i)-(iii) at our disposal, we have that the norm $\|\theta\|_{L^\infty (W^{\frac{k+1}{2},p})}$ is controlled for all $\varepsilon>0$: we have proven spatial regularity. Time regularity follows since we have
$$\frac{ \partial^k}{\partial t^k}\theta(x,t)+\nabla \cdot \left(\frac{ \partial^k}{\partial t^k} (v_\varepsilon\,\theta)\right)(x,t)+\Lambda  \left(\frac{ \partial^k}{\partial t^k} \theta\right)(x,t)=\varepsilon \Delta \left(\frac{ \partial^k}{\partial t^k} \theta\right)(x,t).\qquad\blacksquare$$
\subsection{Maximum principle and Besov regularity}
As a motivation for the Theorem \ref{TheoBesov} below, we rewrite in the following lines the proof of the maximum principle.
\begin{Theoreme}[Maximum Principle]\label{TheoPrincipeMax}
Let $2\leq p<+\infty$ and let $\theta$ be a smooth solution of equation (\ref{SistApprox}). Then we have the following estimation
\begin{equation}\label{PrincipeMax}
\|\theta(\cdot, t)\|_{L^p} \leq \|\theta_0\|_{L^p}.
\end{equation}
\end{Theoreme}
\textit{\textbf{Proof}.} We write
\begin{eqnarray*}
\frac{d}{dt}\|\theta(\cdot, t)\|_{L^p}^p&=&p\int_{\mathbb{R}^n}|\theta|^{p-2}\theta\bigg(\varepsilon \Delta \theta-\nabla \cdot (v_\varepsilon \,\theta)-\Lambda \theta\bigg)dx=p\varepsilon\int_{\mathbb{R}^n}|\theta|^{p-2}\theta\Delta \theta dx+p\int_{\mathbb{R}^n}|\theta|^{p-2}\theta \Lambda \theta dx\\
\end{eqnarray*}
where we used the fact that $div(v)=0$. Thus, we have $\displaystyle{\frac{d}{dt}\|\theta(\cdot, t)\|_{L^p}^p-p\varepsilon\int_{\mathbb{R}^n}|\theta|^{p-2}\theta \Delta \theta dx+p\int_{\mathbb{R}^n}|\theta|^{p-2}\theta \Lambda \theta dx=0},$ and integrating in time we obtain
\begin{equation}\label{Form1}
\|\theta(\cdot, t)\|_{L^p}^p-p\varepsilon\int_{0}^t\int_{\mathbb{R}^n}|\theta|^{p-2}\theta\Delta \theta dxds+p\int_0^t\int_{\mathbb{R}^n}|\theta|^{p-2}\theta \Lambda \theta dxds=\|\theta_0\|_{L^p}^p.
\end{equation}
To finish, we have the next lemma
\begin{Lemme} 
The quantities $\displaystyle{-p\varepsilon\int_{\mathbb{R}^n}|\theta|^{p-2}\theta \Delta \theta dx\qquad \mbox{and} \qquad p\int_0^t\int_{\mathbb{R}^n}|\theta|^{p-2}\theta \Lambda \theta dxds}$ are both positive.
\end{Lemme}
\textit{\textbf{Proof.}} For the first expression, since $e^{\varepsilon s\Delta}$  is a contraction semi-group we have
 $\|e^{\varepsilon s\Delta}f\|_{L^p}\leq \|f\|_{L^p}$ for all $s>0$ and all $f\in L^p(\mathbb{R}^n)$. Thus $F(s)=\|e^{\varepsilon s\Delta}f\|_{L^p}$ is decreasing in $s$; taking the derivative in $s$ and evaluating in $s=0$ we obtain the desired result. For the second expression a proof can be found in \cite{Cordoba} (the positivity lema p.516). However, we will give another proof of this fact with the Theorem \ref{TheoBesov} below.  \hfill$\blacksquare$\\

Getting back to (\ref{Form1}), we have that all these quantities are bounded and positive, so Theorem \ref{TheoPrincipeMax} follows easily. \hfill$\blacksquare$
\begin{Remarque}
\emph{This maximum principle (\ref{PrincipeMax}) is still valid for $1\leq p \leq +\infty$. See \cite{Marchand} for a proof.}
\end{Remarque}
As said in the introduction, the study of expression (\ref{Form1}) above lead us to a result concerning weak solution's regularity which is announced in the Theorem \ref{TheoRegBesov}. More precisely we have 
\begin{Theoreme}[Besov Regularity]\label{TheoBesov}
Let $2\leq p <+\infty$ and let $f:\mathbb{R}^n\longrightarrow \mathbb{R}$ be a function such that
$$\int_{\mathbb{R}^n}|f(x)|^{p-2}f(x)\Lambda f(x)dx<+\infty \quad\mbox{ then  }\quad f\in \dot{B}^{1/p,p}_{p}(\mathbb{R}^n).$$
\end{Theoreme}
\textit{\textbf{Proof}.} To begin with, assume that $f$ is a positive function and let us show that we have the following estimates:
\begin{equation}\label{FormuleEstima1}
\|f\|^p_{ \dot{B}^{2\alpha/p,p}_{p}}\leq C\|f^{p/2}\|^2_{ \dot{B}^{\alpha,2}_{2}}\leq C'\int_{\mathbb{R}^n}|f(x)|^{p-2}f(x)\Lambda f(x)dx
\end{equation}
In this case, we use the following fact:  for $0<\epsilon\leq 1$ and for all $a,b>0$ we have $|a^{\epsilon}-b^{\epsilon}|\leq |a-b|^{\epsilon}$.Hence, applying this fact with $\epsilon=2/p$, $a=f(x)^{p/2}$ and $b=f(y)^{p/2}$, one has $|f(x)-f(y)|\leq  |f(x)^{p/2}-f(y)^{p/2}|^{2/p}$ which implies
$$\|f\|^p_{ \dot{B}^{1/p,p}_{p}}\simeq  \int_{\mathbb{R}^n}\int_{\mathbb{R}^n}\frac{|f(x)-f(y)|^p}{|x-y|^{n+1}}dxdy\leq C \int_{\mathbb{R}^n}\int_{\mathbb{R}^n}\frac{|f(x)^{p/2}-f(y)^{p/2}|^{2}}{|x-y|^{n+1}}dxdy\simeq \|f^{p/2}\|^2_{ \dot{B}^{1/2,2}_{2}}$$
and this give us the first part of (\ref{FormuleEstima1}). For the second part we have
\begin{equation}\label{FormuleEstima2}
\|f^{p/2}\|^2_{ \dot{B}^{1/2,2}_{2}}=\|f^{p/2}\|^2_{ \dot{H}^{1/2}}= \int_{\mathbb{R}^n}|\Lambda^{1/2} f^{p/2}(x)|^2dx=  \int_{\mathbb{R}^n}f^{p/2}(x)\Lambda f^{p/2}(x)dx,
\end{equation}
by the self-adjointness of operator $\Lambda^{1/2}$. We consider now the semi-group $(e^{-\tau \Lambda})_{\tau\geq 0}$. Since $p\geq 2$, using Jensen inequality\footnote{see \cite{PGDCH} for the details concerning the semi-group $(e^{-\tau \Lambda})_{\tau\geq 0}$.} we obtain the estimate $e^{-\tau \Lambda}f\leq \left(e^{-\tau \Lambda}f^{p/2}\right)^{2/p}$. Thus, integrating this inequality we obtain $\|e^{-\tau \Lambda}f\|^p_{L^p}\leq \|e^{-\tau \Lambda}f^{p/2}\|^2_{L^2}$. Finally, taking the derivative with respect to $\tau$ and evaluating in $\tau=0$ one obtains
$$-p\int_{\mathbb{R}^n}|f(x)|^{p-2}f(x)\Lambda f(x)dx \leq -2\int_{\mathbb{R}^n}f^{p/2}(x)\Lambda f^{p/2}(x)dx $$
Hence, getting back to (\ref{FormuleEstima2}) it comes $\|f^{p/2}\|^2_{ \dot{B}^{1/2,2}_{2}}\leq \int_{\mathbb{R}^n}|f(x)|^{p-2}f(x)\Lambda f(x)dx$ and estimates (\ref{FormuleEstima1}) are proven. Let us now prove the general case. For this, we write $f(x)=f_+(x)-f_-(x)$ where $f_\pm(x)$ are positives functions with disjoint support. We have:
\begin{eqnarray}\label{FormDecomp}
\int_{\mathbb{R}^n}|f(x)|^{p-2}f(x)\Lambda f(x)dx&=&\int_{\mathbb{R}^n}f_+(x)^{p-2}f_+(x)\Lambda f_+(x)dx+\int_{\mathbb{R}^n}f_-(x)^{p-2}f_-(x)\Lambda f_-(x)dx\\
&-&\int_{\mathbb{R}^n}f_+(x)^{p-2}f_+(x)\Lambda f_-(x)dx-\int_{\mathbb{R}^n}f_-(x)^{p-2}f_-(x)\Lambda f_+(x)dx<+\infty\nonumber
\end{eqnarray}
We only need to treat the two last integrals, and in fact we just need to study one of them since the other can be treated in a similar way. So, for the third integral we have
\begin{eqnarray*}
\int_{\mathbb{R}^n}f_+(x)^{p-2}f_+(x)\Lambda f_-(x)dx&=&\int_{\mathbb{R}^n}f_+(x)^{p-2}f_+(x)\int_{\mathbb{R}^n}\frac{f_-(x)-f_-(y)}{|x-y|^{n+1}}dydx\\
&=&\int_{\mathbb{R}^n}f_+(x)^{p-2}\int_{\mathbb{R}^n}\frac{f_+(x)f_-(x)-f_+(x)f_-(y)}{|x-y|^{n+1}}dydx
\end{eqnarray*}
However, since $f_+$ and $f_-$ have disjoint supports we obtain the following estimate: 
$$\int_{\mathbb{R}^n}f_+(x)^{p-2}f_+(x)\Lambda f_-(x)dx=-\int_{\mathbb{R}^n}f_+(x)^{p-2}\int_{\mathbb{R}^n}\frac{f_+(x)f_-(y)}{|x-y|^{n+1}}dydx\leq 0$$
This quantity is negative as all the terms inside the integral are positive. With this observation we see that the last terms of (\ref{FormDecomp}) are positive and we have
\begin{eqnarray*}
\int_{\mathbb{R}^n}f_+(x)^{p-2}f_+(x)\Lambda f_+(x)dx+\int_{\mathbb{R}^n}f_-(x)^{p-2}f_-(x)\Lambda f_-(x)dx \leq \int_{\mathbb{R}^n}|f(x)|^{p-2}f(x)\Lambda f(x)dx<+\infty
\end{eqnarray*}
Then, using the first part of the proof we have $f_\pm\in \dot{B}^{1/p,p}_{p}(\mathbb{R}^n)$ and since $f=f_+-f_-$ we conclude that $f$ belongs to the Besov space $\dot{B}^{1/p,p}_{p}(\mathbb{R}^n)$. We have proven the following general estimate
$$\|f\|^p_{ \dot{B}^{1/p,p}_{p}}\leq C\int_{\mathbb{R}^n}|f(x)|^{p-2}f(x)\Lambda f(x)dx\qquad \blacksquare $$
\begin{Remarque}\emph{From this inequality one easily deduces positivity of this last integral. This constitutes another proof for the positivity lemma of \cite{Cordoba} for $2\leq p<+\infty$. Another proof, far more general of the Theorem \ref{TheoBesov} is given in \cite{PGDCH}.}
\end{Remarque}
To obtain weak solutions of (\ref{QG+}) with initial data in $L^p$, $p\geq 2$, we will now pass to the limit by taking $\varepsilon \longrightarrow 0$. We have obtained a family of regular functions $(\theta^{(\varepsilon)})_{\varepsilon >0}\in L^{\infty}([0,T]; L^{p}(\mathbb{R}^n))$ which are solutions of (\ref{SistApprox}) and satisfy the uniform bound 
$$\|\theta^{(\varepsilon)}(\cdot, t)\|_{L^p}\leq \|\theta_0\|_{L^p}$$
Since $L^{\infty}([0,T]; L^{p}(\mathbb{R}^n))= \left(L^{1}([0,T]; L^{q}(\mathbb{R}^n))\right)'$, with $\frac{1}{p}+\frac{1}{q}=1$, we can extract from those solutions $\theta^{(\varepsilon)}$ a subsequence $(\theta_k)_{k\in \mathbb{N}}$ which is $\ast$-weakly convergent to some function $\theta$ in the space $L^{\infty}([0,T]; L^{p}(\mathbb{R}^n))$, which implies convergence in $\mathcal{D}'(\mathbb{R}^+\times\mathbb{R}^n)$. However, this weak convergence is not sufficient to assure the convergence of $(v_\varepsilon\; \theta_k)$ to $v\;\theta$. For this we use the remarks that follows. First, using remark \ref{Rem_Approx} we can consider a sequence $(v_k)_{k\in\mathbb{N}}$ with $v_k$ as in formula (\ref{Forbmoaprox}) such that $v_k \longrightarrow v$ weakly in $bmo$. Secondly, combining (\ref{PrincipeMax}) and Theorem \ref{TheoBesov} we obtain that solutions $\theta_k$ belongs to the space $L^{\infty}([0,T]; L^{p}(\mathbb{R}^n))\cap L^{1}([0,T]; \dot{B}^{1/p,p}_p(\mathbb{R}^n))$ for all $k\in \mathbb{N}$. \\

To finish, fix a function $\varphi\in \mathcal{C}^{\infty}_{0}([0,T]\times \mathbb{R}^n)$. Then we have the fact that $\varphi \theta_k\in  L^{1}([0,T]; \dot{B}^{1/p,p}_p(\mathbb{R}^n))$ and $\partial_t \varphi \theta_k\in  L^{1}([0,T]; \dot{B}^{-N,p}_p(\mathbb{R}^n))$. This implies the local inclusion, in space as well as in time, $\varphi \theta_k\in \dot{W}^{1/p,p}_{t,x}\subset \dot{W}^{1/p,2}_{t,x}$ so we can apply classical results such as the Rellich's theorem to obtain convergence of $v_k\; \theta_k$ to $v\;\theta$. \\

Thus, we obtain existence and uniqueness of weak solutions for the problem (\ref{QG+}) with an initial data in $\theta_0\in L^p(\mathbb{R}^n)$, $2\leq p<+\infty$. Moreover, since such solutions satisfy inequality (\ref{PrincipeMax}) we have that these solutions $\theta(x,t)$ belongs to the space $L^{\infty}([0,T]; L^{p}(\mathbb{R}^n))\cap L^{p}([0,T]; \dot{B}^{1/p,p}_p(\mathbb{R}^n))$. \hfill  $ \blacksquare $ 
\section{Positivity principle}\label{Sect_PrincipeMax}
We prove in this section the Theorem \ref{TheoPrincipePos}. Recall that by hypothesis we have $0\leq \psi_0\leq M$ and $\psi_0\in L^p(\mathbb{R}^n)$ with $1\leq p\leq+\infty$. To begin with, we fix two constants, $\rho, R$ such that $R>2\rho>0$. Then we set $A_{0,R}(x)$ a function equals to $M/2$ over $|x|\leq 2R$ and equals to $\psi_0(x)$ over $|x|>2R$ and we write $B_{0,R}(x)=\psi_0(x)-A_{0,R}(x)$, so by construction we have $$\psi_0(x)=A_{0,R}(x)+B_{0,R}(x)$$ with $\|A_{0,R}(\cdot)\|_{L^\infty} \leq M$ and $\|B_{0,R}(\cdot)\|_{L^\infty} \leq M/2$. Remark that $A_{0,R}, B_{0,R}\in L^p(\mathbb{R}^n)$. Now fix $v\in L^{\infty}([0,T];bmo(\mathbb{R}^n))$ such that $div(v)=0$ and consider the equations
\begin{equation}\label{ForDouble}
\begin{array}{lcl}
\partial_t A_R(x,t)+ \nabla\cdot (v\,A_R)(x,t)+\Lambda A_R(x,t)=0& \qquad\qquad\qquad & \partial_t B_R(x,t)+ \nabla\cdot (v\,B_R)(x,t)+\Lambda B_R(x,t)=0\\
& \begin{Large}and\end{Large} &\\
A_R(x,0)=A_{0,R}(x). &  &B_R(x,0)= B_{0,R}(x).
\end{array}
\end{equation}
Using the maximum principle and by construction we have the following estimates for $t\in [0,T]$:
\begin{eqnarray}
\|A_R(\cdot, t)\|_{L^p}&\leq& \|A_{0,R}\|_{L^p}\leq \|\psi_0\|_{L^p}+CM R^{n/p} \quad (1<p<+\infty)\label{Formula1}\\
\|A_R(\cdot, t)\|_{L^\infty}&\leq & \|A_{0,R}\|_{L^\infty}\leq M.\nonumber\\
\|B_R(\cdot, t)\|_{L^\infty}&\leq & \|B_{0,R}\|_{L^\infty}\leq M/2. \nonumber
\end{eqnarray}
\begin{Lemme}\label{lemmeRecollement} 
The function $\psi(x,t)=A_R(x,t)+B_R(x,t)$,  where $A_R(x,t)$ and $B_R(x,t)$ are solutions of the systems (\ref{ForDouble}), is the unique solution for the problem
\begin{equation}\label{Equation1}
\left\lbrace
\begin{array}{l}
\partial_t \psi(x,t)+ \nabla\cdot (v\,\psi)(x,t)+\Lambda\psi(x,t)=0\\[5mm]
\psi(x,0)=A_{0,R}(x)+B_{0,R}(x).
\end{array}
\right.
\end{equation}
\end{Lemme}
\textit{\textbf{Proof}.} Using hypothesis over $A_R(x,t)$ and $B_R(x,t)$ and the linearity of equation (\ref{Equation1}) we have that the function $\psi_R(x,t)=A_R(x,t)+B_R(x,t)$ is a solution for this equation. Uniqueness is assured by the maximum principle and by the continuous dependence from initial data given in the Corollary \ref{CorDepContinue}, thus we can write $\psi(x,t)$ instead of $\psi_R(x,t)$. \hfill$\blacksquare$\\

To continue, we will need an auxiliary function $\phi \in \mathcal{C}^{\infty}_{0}(\mathbb{R}^n)$ such that $\phi(x)=0$ for $|x|\geq 1$ and $\phi(x)=1$ if $|x|\leq 1/2$ and we set $\varphi(x)=\phi(x/R)$. Now, we will estimate the $L^p$-norm of $\varphi(x)(A_R(x,t)-M/2)$ with $p>n$. 
\begin{Remarque}
\emph{Although some of the following calculations are valid for $1\leq p\leq +\infty$, we will need at the end the fact that $p>n$.}
\end{Remarque}
We write:
\begin{eqnarray}
\partial_t\|\varphi(\cdot)(A_R(\cdot,t)-M/2)\|_{L^p}^p &=& p\int_{\mathbb{R}^n}\big|\varphi(x)(A_R(x,t)-M/2)\big|^{p-2}\big(\varphi(x)(A_R(x,t)-M/2) \big)\nonumber\\[4mm]
& & \times\; \partial_t\big(\varphi(x)(A_R(x,t)-M/2) \big) dx \label{equat1}
\end{eqnarray}
We observe that we have the next identity for the last term above
\begin{eqnarray*}
\partial_t(\varphi(x)(A_R(x,t)-M/2))&=&- \nabla \cdot (\varphi(x) \, v(A_R(x,t)-M/2))-\Lambda(\varphi(x)(A_R(x,t)-M/2))\\[3mm]
&+&(A_R(x,t)-M/2)v\cdot \nabla \varphi(x)+[\Lambda, \varphi]A_R(x,t)-M/2 \Lambda\varphi(x)
\end{eqnarray*}
where we noted $[\Lambda, \varphi]$ the commutator between $\Lambda$ and $\varphi$. Thus, using this identity in (\ref{equat1}) and the fact that $div(v)=0$ we have
\begin{eqnarray}
\partial_t\|\varphi(\cdot)(A_R(\cdot,t)-M/2)\|_{L^p}^p &=&-p\int_{\mathbb{R}^n}\big|\varphi(x)(A_R(x,t)-M/2)\big|^{p-2}\big(\varphi(x)(A_R(x,t)-M/2) \big)\nonumber\\[4mm]
& & \times\; \Lambda(\varphi(x)(A_R(x,t)-M/2))dx\label{equa2}\\[4mm]
&+&p\int_{\mathbb{R}^n}\big|\varphi(x)(A_R(x,t)-M/2)\big|^{p-2}\big(\varphi(x)(A_R(x,t)-M/2) \big)\nonumber\\[4mm]
& & \times\; \left([\Lambda, \varphi]A_R(x,t)-M/2 \Lambda\varphi(x)\right)dx\nonumber
\end{eqnarray}
Remark that integral (\ref{equa2}) is positive so one has
\begin{eqnarray*}
\partial_t\|\varphi(\cdot)(A_R(\cdot,t)-M/2)\|_{L^p}^p &\leq &p\int_{\mathbb{R}^n}\big|\varphi(x)(A_R(x,t)-M/2)\big|^{p-2}\big(\varphi(x)(A_R(x,t)-M/2) \big)\\[4mm]
& & \times\; \left([\Lambda, \varphi]A_R(x,t)-M/2 \Lambda\varphi(x)\right)dx
\end{eqnarray*}
Using Hölder inequality and integrating the previous expression we have
\begin{eqnarray}\label{EquaPrevious}
\|\varphi(\cdot)(A_R(\cdot,t)-M/2)\|_{L^p}^p& \leq &\|\varphi(\cdot)(A_R(\cdot,0)-M/2)\|_{L^p}^p+\int_{0}^t\left\|[\Lambda, \varphi]A_R(\cdot,s)\right\|_{L^p}+ \|M/2 \Lambda\varphi\|_{L^p}ds
\end{eqnarray}
The first term of the right side is null since over the support of $\varphi$ we have identity $A_R(x,0)=M/2$. For the second term $\left\|[\Lambda, \varphi]A_R(\cdot,s)\right\|_{L^p}$ we have the estimate below given by Calder\'on's commutator (see \cite{Grafakos}) and by the maximum principle
\begin{equation*}
\left\|[\Lambda, \varphi]A_R(\cdot,s)\right\|_{L^p}\leq CR^{-1}\|A_R(\cdot,s)\|_{L^p} \leq CR^{-1}\|A_{0,R}\|_{L^p}.
\end{equation*}
Now,  getting back to the last term of (\ref{EquaPrevious}) we have by definition of $\varphi$ the estimate $\|M/2 \Lambda\varphi\|_{L^p}\leq C MR^{n/p-1}$. We thus have
$$\|\varphi(\cdot)(A_R(\cdot,t)-M/2)\|_{L^p}^p\leq CR^{-1}\int_{0}^t\bigg(\|A_{0,R}\|_{L^p}+MR^{n/p}\bigg)ds.$$
Observe that we have at our disposal estimate (\ref{Formula1}), so we can write
$$\|\varphi(\cdot)(A_R(\cdot,t)-M/2)\|_{L^p}^p\leq Ct R^{-1}\left(\|\psi_0\|_{L^p}+MR^{n/p}\right)$$
Using again the definition of $\varphi$ one has
$$\left(\int_{B(0,\rho)}|A_R(\cdot,t)-M/2|^{p}dx\right)\leq Ct R^{-1}\left(\|\psi_0\|_{L^p}+MR^{n/p}\right).$$
Thus, if $R\longrightarrow +\infty$ and since $p>n$, we have $\alpha(x,t)=M/2$ over $B(0,\rho)$.\\ 

Hence, by construction we have  $\psi(x,t)=A_R(x,t)+B_R(x,t)$ where $\psi$ is a solution of $(T)_{\alpha}$ with initial data $\psi_0=A_{0,R}+B_{0,R}$, but, since over $B(0,\rho)$ we have $A(x,t)=M/2$ and $\|B(\cdot,t)\|_{L^\infty}\leq M/2$, one finally has the desired estimate $0\leq \psi(x,t)\leq M$. \hfill$\blacksquare$
\section{Existence of solutions with a $L^\infty$ initial data}\label{SecLinfty}
The proof given before for the maximum principle allows us to obtain the existence of solutions for fractional diffusion transport equation (\ref{QG}) when the initial data $\theta_0$ belongs to the space $L^\infty(\mathbb{R}^n)$. Indeed, let us fix $\theta_0^R=\theta_0 \mathds{1}_{B(0,R)}$ with $R>0$ so we have $\theta_0^R\in L^p(\mathbb{R}^n)$ for all $1\leq p\leq +\infty$. Following section \ref{SecExiUnic}, there is a unique solution $\theta^R$ for the problem
\begin{equation*}
\left\lbrace
\begin{array}{l}
\partial_t \theta^R+ \nabla\cdot (v\theta^R)+\Lambda\theta^R=0\\[5mm]
\theta^R(x,0)=\theta_0^R(x)\\[5mm]
div(v)=0 \quad \mbox{ and } v\in L^{\infty}([0,T];  bmo(\mathbb{R}^n)).
\end{array}
\right.
\end{equation*}
such that $\theta^R\in L^\infty([0,T]; L^p(\mathbb{R}^n))$. By the maximum principle we have $\|\theta^R(\cdot, t)\|_{L^p}\leq\|\theta^R_0\|_{L^p}\leq v_n\|\theta_0\|_{L^\infty}R^{n/p}$. Taking the limit $p\longrightarrow +\infty$ and making $R\longrightarrow +\infty$ we finally have $\|\theta(\cdot, t)\|_{L^\infty}\leq C\|\theta_0\|_{L^\infty}$. This shows that for an initial data $\theta_0\in  L^\infty(\mathbb{R}^n)$ there exists an associated solution $\theta\in L^\infty([0,T];L^\infty(\mathbb{R}^n))$.
\section{Application to the 2D-quasi-geostrophic equation}\label{SecAppli}
We have worked so far with a velocity given by a general function $v\in L^{\infty}([0,T];bmo(\mathbb{R}^n))$, let us now treat super-critical case of the 2D-quasi-geostrophic equation with $u=(-R_2\theta,R_1\theta)$; where $R_j$ are the Riesz transforms. Fix $\theta_0$ an initial data belonging to $L^p\cap L^\infty(\mathbb{R}^2)$, with $p\geq 2$. Following \cite{Marchand} we have the existence of a solution $\theta(\cdot, t)$ for the equation $(QG)_{1/2}$ with $\theta(\cdot,t)\in L^p\cap L^\infty(\mathbb{R}^2)$ for $t\in [0,T]$. Since Riesz transforms are bounded in $L^p$ and since they are bounded from $L^\infty$ into $BMO$, we have a uniform bound of the velocity $u$ in terms of the $bmo$ norm: we can apply the Theorem \ref{Theo1} to obtain Hölder regularity for the solution of 2D-quasi-geostrophic equation.
\pagebreak

\quad\\[5mm]

\begin{flushright}
\begin{minipage}[r]{80mm}
Diego \textsc{Chamorro}\\[3mm]
Laboratoire d'Analyse et de Probabilités\\ 
Université d'Evry Val d'Essonne\\[2mm]
91025 Evry Cedex\\[2mm]
diego.chamorro@univ-evry.fr
\end{minipage}
\end{flushright}

\end{document}